\documentclass{amsart}
\usepackage{amsmath, amssymb, stmaryrd, amsthm, mathtools}
\usepackage[all,knot]{xy}

\usepackage{tikz}
\usetikzlibrary{%
  shapes,
  decorations.shapes,
  decorations.fractals,
  decorations.markings,
  shadows
}
\usetikzlibrary{arrows}
\definecolor{maroon}{HTML}{BB0000}

\usetikzlibrary{decorations.markings,backgrounds,calc}

\tikzset{middlearrow/.style={
        decoration={markings,
            mark= at position 0.5 with {\arrow[scale=1.25]{#1}} ,
        },
        postaction={decorate}
    }
}

\usepackage{url,verbatim,hyperref,geometry,parskip}

\newcommand{\mcB}{\mathcal{B}}
\newcommand{\mcC}{\mathcal{C}}
\newcommand{\mcD}{\mathcal{D}}
\newcommand{\mcE}{\mathcal{E}}

\newcommand{\mcM}{\mathcal{M}}

\newcommand{\mcQ}{\mathcal{Q}}

\newcommand{\mcS}{\mathcal{S}}

\newcommand{\s}{\sigma}

\newcommand{\w}{\omega}

\newcommand{\CC}{\mathcal{C}}

\newcommand{\Tt}{U(1)}

\renewcommand{\bar}[1]{\overline{#1}}

\newcommand{\0}{\mathbf{0}}
\renewcommand{\1}{\mathbf{1}}

\DeclareMathOperator{\Hom}{Hom}

\DeclareMathOperator{\Id}{Id}

\DeclareMathOperator{\Vect}{Vect}
\renewcommand{\Vec}{\text{Vec}}
\DeclareMathOperator{\ord}{ord}

\DeclareMathOperator{\Rep}{Rep}

\DeclareMathOperator{\sVec}{sVec}

\DeclareMathOperator{\ii}{i}

\numberwithin{equation}{section}
\newtheorem{theorem}{Theorem}[section]

\newtheorem{corollary}[theorem]{Corollary}
\newtheorem{lemma}[theorem]{Lemma}

\newtheorem{question}[theorem]{Question}

\newtheorem{conjecture}[theorem]{Conjecture}

\newtheorem{prop}[theorem]{Proposition}
\newtheorem{proposition}[theorem]{Proposition}

\theoremstyle{definition}
\newtheorem{definition}{Definition}[section]

\theoremstyle{definition}
\newtheorem{example}{Example}[section]

\theoremstyle{remark}
\newtheorem{remark}{Remark}[section]

\newcommand{\II}{\mathcal{I}}

\newcommand{\ot}{{\otimes}}

\newcommand{\one}{\mathbf{1}}
\newcommand{\triv}{\mathbf{1}}

\newcommand{\Z}{\mathbb{Z}}

\newcommand{\mC}{\mathcal{C}}

\newcommand{\EE}{\mathcal{E}}

\newcommand{\mcd}{\mathcal{D}}

\newcommand{\mS}{\mathcal{S}}

\newcommand{\qcal}{\mathcal{Q}}

\newcounter{commentcounter}

\renewcommand\o{\otimes}
\newcommand{\BF}{\mathbb{F}}
\newcommand{\mG}{\mathcal{G}}
\newcommand{\mD}{\mathcal{D}}
\newcommand{\id}{\operatorname{Id}}
\renewcommand{\Rep}{\operatorname{Rep}}

\begin{document}

\title{Fermionic modular categories and the 16-fold Way}
\date{\today}
\author[Bruillard]{Paul Bruillard}
\email{pjb2357@gmail.com}
\address{Pacific Northwest National Laboratory, 902 Battelle Boulevard,
Richland, WA U.S.A.}

\author[Galindo]{C\'{e}sar Galindo}
\email{cn.galindo1116@uniandes.edu.co}
\address{Departamento de Matem\'aticas, Universidad de los Andes, Bogot\'a, Colombia.}

\author[Hagge]{Tobias Hagge}
\email{tobias.hagge@pnnl.gov}
\address{Pacific Northwest National Laboratory, 902 Battelle Boulevard,
Richland, WA U.S.A.}

\author[Ng]{Siu-Hung Ng}
\email{rng@math.lsu.edu}
\address{Department of Mathematics, Louisiana State University, Baton Rouge, LA
    U.S.A.}
\author[Plavnik]{Julia Yael Plavnik}
\email{julia@math.tamu.edu}
\address{Department of Mathematics,
    Texas A\&M University,
    College Station, TX
    U.S.A.}
\author[Rowell]{Eric C. Rowell}
\email{rowell@math.tamu.edu}
\address{Department of Mathematics,
    Texas A\&M University,
    College Station, TX
    U.S.A.}
\author[Wang]{Zhenghan Wang}
\email{zhenghwa@microsoft.com}
\address{Microsoft Research Station Q and Department of Mathematics,
    University of California,
    Santa Barbara, CA
    U.S.A.}

\begin{abstract}
We study spin and super-modular categories systematically as inspired by fermionic topological phases of matter, which are always fermion parity enriched and modelled by spin TQFTs at low energy.  We formulate a $16$-fold way conjecture for the minimal modular extensions of super-modular categories to spin modular categories, which is a categorical formulation of gauging the fermion parity.  We investigate general properties of super-modular categories such as fermions in twisted Drinfeld doubles, Verlinde formulas for naive quotients, and explicit extensions of $PSU(2)_{4m+2}$ with an eye towards a classification of the low-rank cases.

\end{abstract}
\thanks{The results in this paper were mostly obtained while all authors except the third were at the American Institute of Mathematics during August 10-14,2015, participating in a SQuaRE.  We would like to thank AIM for their hospitality and encouragement. C. Galindo was partially supported by the FAPA funds from vicerrectoria de investigaciones de la Universidad de los Andes, S.-H. Ng  by NSF grant DMS-1501179, J.Plavnik by NSF grant DMS-1410144,  CONICET, ANPCyT, and Secyt-UNC, E. Rowell by NSF grants DMS-1108725 and DMS-1410144, and Z. Wang by NSF grants DMS-1108736 and 1411212. Z.W. thanks N. Read for an earlier collaboration on a related topic \cite{RW} in which some of the materials in this paper were discussed.}
\maketitle

\section{Introduction}

The most important class of topological phases of matter is two dimensional electron liquids which exhibit the fractional quantum Hall effect (see \cite{DFNSS} and references therein).  Usually fractional quantum Hall liquids are modelled by Witten-Chern-Simons topological quantum field theories (TQFTs) at low energy based on bosonization such as flux attachment.  But subtle effects due to the fermionic nature of electrons are better modelled by refined theories of TQFTs (or unitary modular categories) such as spin TQFTs (or fermionic modular categories) \cite{BW,RW,GWW}.  
In this paper, we study a refinement of unitary modular categories to spin modular categories \cite{Bl,RW} and their local  sectors---super-modular categories \cite{Bon,DGNO1,LW,Sawin}.

Let $f$ denote a fermion in a fermionic topological phase of matter, and $\mathbf{1}$ be the ground state of an even number of fermions.  Then in fermion systems like the fractional quantum Hall liquids, $f$ cannot be distinguished topologically from $\mathbf{1}$ as anyons, so in the low energy effective theory we would have $f\cong \mathbf{1}$.  We would refer to this mathematical identification $f\cong \mathbf{1}$ as the condensation of fermions.  This line of thinking leads to a mathematical model as follows:  the local sector of a fermionic topological phase of matter will be modelled by a super-modular category $\mcB$---a unitary pre-modular category such that every non-trivial transparent simple object is isomorphic to the fermion $f$.   To add the twisted or defect sector associated to fermion parity, we will extend the super-modular category $\mcB$ to a unitary modular category $\mcC$ with the smallest possible dimension $D_\mcC^2=2D^2_\mcB$. Such a unitary modular category has a distinguished fermion $f$ and will be called a spin modular category.  We will also say that $\mcC$ covers the super-modular category $\mcB$.  If the fermion $f$ in $\mcC$ is condensed, then we obtain a ferminonic quotient $\mcQ$ of $\mcC$.  But an abstract theory of such fermionic modular categories $\mcQ$ has not been developed.  Given a super-modular category $\mcB$, it is open whether or not there will always be a covering spin modular category.  If a covering theory exists, then it is not unique.  One physical implication is that a super-modular category alone is not enough to characterize a fermionic topological order, which is always fermion parity enriched.  We need the full spin modular category to classify fermionic topological orders such as fermionic fractional quantum Hall states \cite{RW}.  In this paper, we study the lifting of super-modular categories to their spin covers.  

Fermion systems have a fermion number operator $(-1)^F$ which leads to the fermion parity: eigenstates of $(-1)^F$ with eigenvalue $+1$ are states with an even number of fermions and  eigenstates of $(-1)^F$ with eigenvalue $-1$ are states with an odd number of fermions.  This fermion parity is like a $\mathbb{Z}_2$-symmetry in many ways, but it is not strictly a symmetry because fermion parity cannot be broken.  Nevertheless, we can consider the gauging of the fermion parity (compare with \cite{BBCW,CGPW}).  In our model, the gaugings of the fermion parity are the minimal extensions of the super-modular category $\mcB$ to its covering spin modular categories $\mcC$.  We conjecture that a minimal modular extension always exists, and there are exactly $16$ such minimal extensions of super-modular categories.  We will refer to this conjecture as the $16$-fold way conjecture \ref{16foldway}.  We prove that if there is one minimal extension, then there are exactly $16$ up to Witt equivalence.  A stronger result \cite[Theorem 5.3]{LKW} replaces Witt equivalence by ribbon equivalence.  Therefore, the difficulty in resolving the $16$-fold way conjecture lies in the existence of at least one minimal extension.  We analyze explicitly the minimal modular extensions of the super-modular categories $PSU(2)_{4m+2}, m\geq 0$ using a new construction called zesting.  Zesting applies to more general settings and is our main technical contribution.  Given a modular closure using zesting we can constructs eight
new closures each one with different central charge.

The contents of the paper are as follows.  In section 2, we discuss basic properties of spin modular categories, and describe explicitly fermions in symmetric fusion categories and twisted Drinfeld doubles.  In section 3, we formulate the $16$-fold way conjecture.  We provide support for the conjecture by proving the $16$-fold way for Witt classes given existence, and analyzing explicitly the $16$-fold way for $PSU(2)_{4m+2}, m\geq 0$. Finally, in section 4, we discuss spin TQFTs.

\section{Spin modular categories}

We will work with unitary categories over the complex numbers $\mathbb{C}$ in this paper due to our application to topological phases of matter.  Many results can be generalized easily to the non-unitary setting and ground fields other than $\mathbb{C}$.  Spin modular categories without unitarity were first studied in \cite{Bl}.

\subsection{Fermions}

Let ${\mcB}$ be a unitary ribbon fusion category (URFC), and $\Pi_\mcB$ the set of isomorphism classes of simple objects of $\mcB$, called the {\emph{label set}}.  A URFC is also called a unitary pre-modular category or a unitary braided fusion category.  Given a label $\alpha \in \Pi_\mcB$, we will use $X_\alpha$ to denote a representative object with label $\alpha$.  In general, it is important to distinguish between labels and the representative simple objects in their classes.  But sometimes, we will use $\alpha$ for both the label and a simple object in the class $\alpha$.  A chosen unit of $\mcB$ will be denoted by $\one$, and its label by $0$.  Tensor product $\otimes$ of objects will sometimes be written simply as multiplication.

Given a URFC $\mcB$, let $d_\alpha=\dim(X_\alpha)$ and $\theta_\alpha$ be the quantum dimension and twist of the label $\alpha$, respectively.  The entries of the unnormalized $S$-matrix will be $\tilde{s}_{ij}$, and the normalized $S$-matrix is $s=\frac{\tilde{s}}{{D}}$, where ${D}^2=\dim(\mcB)=\sum_{\alpha\in\Pi_\mcB }d_\alpha^2$.  Braiding of two objects $X,Y$ will be denoted by $c_{X,Y}$.  When $XY$ is simple, then $c_{X,Y}\cdot c_{Y,X}$ is $\lambda_{XY}\cdot \Id_{XY}$ for some scalar $\lambda_{XY}$.  If $X_i, X_j$ and $X_iX_j$ are all simple, then $\lambda_{ij}=\frac{\tilde{s}_{ij}}{d_id_j}$.

\begin{definition}

\begin{enumerate}

\item A \textit{fermion} in a URFC is a simple object $f$ such that $f^2=\one$ and $\theta_f=-1$.

\item A \textit{spin modular category} is a pair $(\mC, f)$, where $\mC$ is a unitary modular category (UMC), and $f$ is a fixed fermion.
\end{enumerate}

\end{definition}

\begin{remark}\label{Remark definition of fermion}
If $X$ is an invertible object in a URFC $\mcB$, then $c_{X,X} = \theta_{X}\id_{X \o X}$, see, for example, \cite[Appendix E.3]{kitaev}. An equivalent definition of a fermion in a URFC $\mcB$ is an object $f$ such that $f^2=\one$ and $c_{f,f}=-1$. Note this definition makes sense in an arbitrary unitary braided fusion category.
\end{remark}

\subsection{Fermions in unitary symmetric fusion categories and twisted Drinfeld doubles}

Recall that a braided fusion category $(\mcC,c)$ is called symmetric if $c_{Y,X}c_{X,Y}=\id_{X\otimes Y}$ for all $X,Y\in\mcC$.  

The fusion category $\operatorname{Rep}(G)$ of complex finite dimensional representations of a finite group $G$ with the canonical braiding $c_{X,Y}(x\otimes y)=y\otimes x$ is an example of a symmetric tensor category called a \textit{Tannakian} fusion category.  More general symmetric fusion categories are constructed as the category of representations of a finite super-group. A finite super-group is a pair $(G,z)$, where $G$ is a finite group and $z$ is a central element of order$\leq 2$. An irreducible representation of $G$ is odd if $z$ acts as the scalar $-1$, and is even if $z$ acts as the identity. If the degree of a simple object $X$ is denoted by $|X| \in \{0,1\}$, then the braiding of two simple objects $X$, $Y$ is
$$c'_{X,Y}(x\otimes y)= (-1)^{|X||Y|}y\otimes x.$$The category $\operatorname{Rep}(G)$ with the braiding $c'$ is called a \textit{super-Tannakian} category, and denoted by $\Rep(G,z)$. 
Any (pseudo-)unitary fusion category has a unique pivotal spherical structure so that $d_\alpha>0$ for all simple objects $\alpha$. 
With respect to this choice we have $\theta_V=-\id_V$ for any odd simple $V \in \Rep(G,z)$, so that $\Rep(G,z)$ is Tannakian exactly when $z=1$.  By \cite[Corollaire 0.8]{De}, every symmetric fusion category is equivalent to a super-Tannakian (possibly Tannakian) category.  A key example of a unitary super-Tannakian category is $\sVec$, the category of super-vector spaces, which has $\tilde{s}_{\sVec}=\left(\begin{smallmatrix} 1 & 1\\1&1\end{smallmatrix}\right)$ and $T_{\sVec}=\left(\begin{smallmatrix} 1 &0 \\0&-1\end{smallmatrix}\right)$.
\begin{remark}
In the literature $\sVec$ usually refers to the symmetric fusion category.  There are two possible pivotal spherical structures that render $\sVec$ a symmetric ribbon category: one gives the unitary version we study, the other has trivial twists but the non-trivial simple object has dimension $-1$.  For us, $\sVec$ will always be the unitary symmetric ribbon category.
\end{remark}
\begin{proposition}\label{prop: symmetric contains fermion}
A symmetric fusion category $\mathcal{C}$ admits a fermion if and only if it is of the form $\operatorname{Rep}(G)\boxtimes \sVec$.
\end{proposition}

\begin{proof}

By Remark \ref{Remark definition of fermion}, Tannakian categories do not admit fermions.
Fermions in a super-Tannakian category are in one-to-one correspondence with group homomorphisms  $\chi:G\to \{1,-1\}$ such that $\chi(z)=-1$. 
Thus, if a super group $(G,z)$ admits a fermion, then $G\cong G/\langle z\rangle\times \mathbb{Z}/2\mathbb{Z}$. It follows that $\Rep(G/\langle z\rangle)\boxtimes \sVec\cong \Rep(G,z)$ as symmetric fusion categories.
\end{proof}

\begin{remark} There are  unitary non-Tannakian symmetric categories that do not admit a fermion, i.e. not of the form $\operatorname{Rep}(G)\boxtimes\sVec$.
One example is the super-Tannakian category $\Rep(\mathbb{Z}_4, 2)$, in which there is a pair of dual simple objects with twist $\theta=-1$, while the other non-trivial object is a boson.
\end{remark}
Let $G$ be a finite group and $w \in Z^3(G, \Tt)$. Define 

and  
\begin{eqnarray}
\label{beta-defn}
\beta_a(x,y)&=& \frac{w(a,x,y)w(x,y,y^{-1}x^{-1}
axy)}{w(x,x^{-1}ax,y)},\\
\label{gamma-defn}
\gamma_a(x,y) &=& \frac{w(x,y,a)w(a,a^{-1}xa,a^{-1}ya)}
{w(x,a,a^{-1}ya)},
\end{eqnarray}
for all $a,x,y\in G$. 
Since $w$ is a 3-cocycle, we have
\begin{equation}
\label{beta-reln}
\beta_a(x,y)\beta_a(xy,z) = \beta_a(x,yz)\beta_{x^{-1}ax}(y,z)
\end{equation}
for all $a, x,y,z\in G$. Therefore, for any $a\in G$ the restriction $\beta_a|_{C_G(a)}$ is a 2-cocycle. 

%
Let us recall the description of the UMC $\operatorname{Rep}(D^w(G))$---the category of representations of the twisted Drinfeld double defined by Dijkgraaf, Pasquier and Roche in
\cite[Section 3.2]{Dijkgraaf}.

An object is a $G$-graded finite dimensional Hilbert space $\mathcal{H}=
\bigoplus_{k\in G}\mathcal H_k$ and a twisted $G$-action, $\rhd : G \to U(\mathcal H)$ such that

\begin{itemize}
  \item $\sigma \rhd \mathcal H_k = \mathcal H_{\sigma \cdot k}$
  \item $\sigma \rhd (\tau \rhd h_k) = \beta_k(\sigma,\tau)(\sigma \tau) \rhd h_k$
  \item $e\rhd h = h$
\end{itemize}
for all $\sigma, \tau, k\in G, h_k\in \mathcal H_k$. Morphisms in the category are linear maps that preserve  the grading and the twisted action, i.e., a linear map  $f:\mathcal H\to \mathcal H'$ is a morphism if
\begin{itemize}
  \item $f(\mathcal H_k)\subset \mathcal H'_k$,
  \item $f(\sigma \rhd h)= \sigma \rhd f(h)$
\end{itemize}for all $\sigma,k \in G$ and $h\in \mathcal H$.

The monoidal structure on $\operatorname{Rep}(D^w(G))$ is defined as follows: let $\mathcal H$ and $\mathcal H'$ be objects in $\operatorname{Rep}(D^w(G))$, then the tensor product of Hilbert spaces $\mathcal H\otimes \mathcal H'$ is an object in $\operatorname{Rep}(D^w(G))$ with $G$-grading $(\mathcal H\otimes \mathcal H')_k =\bigoplus_{x,y \in G: xy = k}\mathcal H_x\otimes \mathcal H'_y$ and twisted $G$-action
$$\sigma\rhd (h_x\otimes h_y'):= \gamma_\sigma(x,y) (\sigma\rhd h_x\otimes \sigma \rhd h_y'),$$
for all $ \sigma, x ,y \in G,$ $h_x \in \mathcal H_x$ and $ h_y' \in \mathcal H_y'$.

Now, for $\mathcal H$, $\mathcal H'$ and $\mathcal H''$ objects
in $\operatorname{Rep}(D^w(G))$ the associativity constraint  $$\Theta:
(\mathcal H\otimes \mathcal H')\otimes \mathcal H''\to \mathcal
H\otimes (\mathcal H'\otimes \mathcal H''),$$ for the monoidal
structure $\otimes$ is defined by $$\Theta((h_x\otimes
h_y')\otimes h_z'')= w(x,y,z)h_x\otimes (h_y'\otimes
h_z'')$$ for all $x, ,y, z\in G$, $h_x\in \mathcal H_x$, $h_y'\in
\mathcal H_y'$ and $h_z''\in \mathcal H_z''$.

The unit object $\underline{\mathbb{C}}$ is defined as the
one dimensional Hilbert space $\mathbb C$ graded only at the unit
element $e\in G$, endowed with trivial $G$-action.

Finally, for  $\mathcal H$ and $\mathcal H'$  objects
in $\operatorname{Rep}(D^w(G))$, the braiding is defined by $$c_{\mathcal H,\mathcal H'}(h_x\otimes h_{y})=x\rhd h_y\otimes h_x,$$
for all $x, y\in G$, $h_x\in \mathcal H_x$ and $h_y'\in
\mathcal H_y'$.

The invertible objects in $\Rep(D^w(G))$ can be parametrized as follows.
If $z\in Z(G)$ then $\beta_z(-,-)\in Z^2(G,\Tt)$.
 Define $Z_w(G)$ as those $z\in Z(G)$ such that $\beta_z(-,-)$ is a 2-coboundary.  Then there exists an 
$\eta:G\to \Tt$ such that $\frac{\eta(\sigma)\eta(\tau)}{\eta(\sigma\tau)}=\beta_z(\sigma,\tau)$ for all $\sigma,\tau \in G$. There is a correspondence between invertible objects in $\Rep(D^w(G))$ and pairs $(\eta,z)$ where $z\in Z_w(G)$ and $\eta$ is as above.  The tensor product is given by $(\eta,z)\otimes (\eta',z')=(\gamma_{(-)}(z,z')\eta\eta',zz')$. In fact, by \cite[Prop.5.3]{MasNg}, the group $\mS$ of invertible objects of $\Rep(D^w(G))$ fits into the exact sequence 
$$
1 \to \hat{G} \to \mS \to Z_w(G)  \to 1
$$
where $\hat{G}$ is the group of linear characters of $G$.   
\begin{proposition}\label{fermions in D(G)}
There is a correspondence between fermions in $\operatorname{Rep}(D^w(G))$ and pairs $(\eta,z)$, where $\eta: G\to U(1)$ and
\begin{itemize}
    \item[(a)] $z\in Z(G)$ of order two,
    \item[(b)] $\frac{\eta(\sigma)\eta(\tau)}{\eta(\sigma\tau)}=\beta_z(\sigma,\tau)$ for all $\sigma,\tau \in G$,
    \item[(c)] $\gamma_z(x,z)\eta(x)^2=1$ for all $x\in G$.
    \item[(d)] $\eta(z)=-1$.
\end{itemize}
\end{proposition}
\begin{proof}
It follows easily from the definition of $\operatorname{Rep}(D^w(G))$.
\end{proof}

If $w=1$, then fermions in $\operatorname{Rep}(D(G))$ correspond just with pairs $(\chi,z)$, where $\chi:G\to \{1,-1\}$ and $z\in G$ a central element of order two such that $\chi(z)=-1$. Then as in Proposition \ref{prop: symmetric contains fermion}, $G\cong \overline{G} \times \langle z\rangle$, where $\overline{G}:= G/\langle z\rangle$ and $\operatorname{Rep}(D(G))\cong \operatorname{Rep}(D(\overline{G}))\boxtimes \operatorname{Rep}(D(\mathbb{Z}_2))$.

If $w$ is not a coboundary and $z\in G$ is a central element of order two, we would like to know is there is $\eta :G\to \Tt$ such that $(\eta,z)$ is a fermion. 

By (d) and (c) of Proposition \ref{fermions in D(G)} if $\operatorname{Rep}(D^{w}(G))$ has a fermion, then $w(z,z,z)=1$. Thus, the first obstruction is that  $w(z,z,z)=1$ or equivalently that the restriction of $w$ to $\langle z\rangle$ is trivial. 

The second obstruction is that the cohomology class of $\beta_z(-,-)\in Z^2(G,\Tt)$ vanishes. Let $\eta:G\to \Tt$ such that $\delta_G(\eta)= \beta_z(-,-)$, then $(\eta,z)$ represents an invertible object in $\operatorname{Rep}(D^w(G))$ and $\beta_z(-,z)\eta(-)^2:G\to \Tt$ is a linear character. The character $\beta_z(-,z)\eta(-)^2$ can be seen as an element in $Z^2(\mathbb{Z}_2,\operatorname{Hom}(G,\Tt)).$ Its cohomology class is zero if and only if there is a linear character $\mu: G\to \Tt$ such that $\mu^2=\beta_z(-,z)\eta(-)^2$ and in this case  $(\eta\mu^{-1},z)$ defines a invertible object in $\operatorname{Rep}(D^{w}(G))$ of order two. 

Now, if $(\eta,z)$ and $(\eta',z)$  are two invertible objects of order two, $\beta_z(-,z)\eta\eta':G\to \{1,-1\}$ is a bicharacter, that is, the set of equivalence classes of invertible objects of order two of the form $(\eta,z)$ with $z$ fixed, is a torsor over $\operatorname{Hom}(G,\{1,-1\})$.  

Finally, (recall that  $w(z,z,z)=1$) if $(\eta,z)$ is an  invertible object of order two, then $\eta(z)\in \{1,-1\}$. If $\eta(z)=1$, the pair $(\eta,z)$ defines a boson and if there exists $\chi:G\to \{1,-1\}$ with $\chi(z)=-1$ the pair $(\chi\eta,z)$ is a fermion.

\begin{example}
Let $G$ be the finite group $\operatorname{SL}(2, \BF_5)$. Then the center $Z(G) = \{\pm I\}$ and we have the exact sequence
  $$
   1 \to Z(G) \to G \to \operatorname{PSL}(2, \BF_5) \to 1\,.
  $$
  Note that $\operatorname{PSL}(2, \BF_5)$ is isomorphic to the simple group $A_5$. Then $\Rep(G, z)$ is super-Tannakian, and the even part of $\Rep(G, z)$ is equivalent to $\Rep(A_5)$ as braided monoidal categories. Since $G$ is a perfect group, $\Rep(G)$ has no linear characters and hence $\Rep(G, z)$ has no fermions. Moreover, every simple object of $\Rep(G,z)$ is self-dual.
  
  By \cite[Prop. 5.2]{MasNg}, for any $w \in Z^3(G, U(1))$, the group of invertible objects of the modular category $\mC=\Rep(D^w(G))$ is isomorphic to $\Z_2$. In particular, $\mC$ has a unique nontrivial invertible object $X$, and the subcategory $\mG$,  generated by the invertible simple objects of  $\mC$, is equivalent to $\Vect(Z(G), w)$ as fusion categories. Then, by \cite[Thm. 5.5]{MasNg},  the ribbon subcategory $\mG$ is modular if and only if the restriction of $w$ on $Z(G)$ is not a coboundary. Since $G$ is the binary icosahedral group, $H^3(G, U(1)) \cong \Z_{120}$. In particular, $G$ is a periodic group (cf. \cite[Chap. XII, 11]{CE}), and so the restriction map $\mbox{res}:H^3(G, U(1)) \rightarrow H^3(Z(G), U(1))$ is surjective. Since $Z(G) \cong \Z_2$, $\mbox{res}(\w'^2)$ is trivial  for any $\w'\in H^3(G, U(1))$.   Therefore, the restriction of $w$ on $Z(G)$ is a coboundary if and only if the order of the cohomology class $\w$ of $w$ in $H^3(G, U(1))$ is a multiple of $8$. 
  
  Suppose $w$ is a representative of $\w \in H^3(G, U(1))$ with $8 \nmid \ord(\w)$, and let $\mD$ be the centralizer of $\mG$ in  $\mC$. Then, by \cite[Thm. 5.5]{MasNg}, $\mG \subseteq \mD$.  If $4 \mid \ord(\w)$, then $\mG$ is equivalent to $\sVec$ (cf. \cite[p243]{MasNg}) and $\mD$ is a super-modular category. Moreover, $\mC$ is the modular closure of $\mD$. If $4\nmid \ord(\w)$, then $\mG$ is Tannakian.
\end{example}

\begin{example}
Let $G$ be a non-abelian group of order eight (dihedral or quaternions) and $z\in Z(G)$ the non trivial central element. By Proposition \ref{prop: symmetric contains fermion}, the symmetric category $\operatorname{Rep}(G,z)$ does not have fermions. Note that the two-dimensional simple representation of $G$ is a self-dual object with twist $\theta=-1$

Let $C_0$ be a cyclic subgroup of $H^3(G, U(1))$ of maximal order $n$. Then $H^3(G, U(1))=C_0\oplus C_1$ for some subgroup $C_1$ of $H^3(G, U(1))$. In fact, $n=4$ if $G$ is the dihedral group and $n=8$ if $G$ is the quaternion group.

Similar to the preceding example, whether or not $\Rep(D^w(G))$ admits a ribbon subcategory equivalent to the semion or $\sVec$ is determined by the order of the coset $\w C_1$ in $H^3(G, U(1))/C_1$, where $\w \in H^3(G, U(1))$ denotes the cohomology class of $w$. The modular category $\Rep(D^w(G))$ admits a semion modular subcategory if and only if $\ord(\w C_1) = n$. The super vector space $\sVec$ is a ribbon subcategory of $\Rep(D^w(G))$ if and only if $\ord(\w C_1) = n/2$ (cf. \cite[Tbl. 2]{GoMN}). Since the group of invertible objects is isomorphic to $\Z_2^3$, if $\Rep(D^w(G))$ admits a ribbon subcategory   equivalent to the semion or $\sVec$, there are exactly four such subcategories.
\end{example}

\begin{example}
Let $w'\in Z^3(\mathbb{Z}_2,U(1))$, given by $w'(1,1,1)=-1$ and $\pi:\mathbb{Z}_4\to \mathbb{Z}_2$ the non-trivial epimorphism. Define $w'\in Z^3(\mathbb{Z}_4,U(1))$ by $w=\pi^*(w')$. If we define $\eta_{\pm}:\mathbb{Z}_4\to U(1), \eta_{\pm}(1)=\eta_{\pm}(3)= \pm i, \eta(2)= -1$, the pairs $(\eta_{\pm},2)$ define two fermions in $\Rep(D^w(\mathbb{Z}_4))$. Note that unlike the case $w=1$, the existence of a fermion over $z$ does not imply that the exact sequence $0\to \langle z\rangle \to G\to Q\to 1$ splits. 
\end{example}

\subsection{General Properties}
The presence of a fermion in a URFC implies several useful properties.
\begin{proposition}\label{fermionproperties}

Let $f$ be a fermion in a URFC $\mcB$, then

\begin{enumerate}

\item Tensoring with $f$ induces an action of $\mathbb{Z}_2$ on the equivalence classes of simple objects.

\item For any label $\alpha$, $\tilde{s}_{f,\alpha}=\epsilon_\alpha d_\alpha$, where $\epsilon_\alpha=\pm 1$ (equivalently, $c_{f,\alpha}c_{\alpha,f}=\epsilon_\alpha \Id_{\alpha\ot f}$). Moreover, $\epsilon_f=1$.

\item $\theta_{f \alpha}=-\epsilon_\alpha \theta_\alpha$.\label{thetasign}

\item $\tilde{s}_{f\alpha,j}=\epsilon_j \tilde{s}_{\alpha,j}$.\label{smatrixsign}

\item $\epsilon_{f\alpha}=\epsilon_\alpha$.  In general, if the fusion coefficient $N_{ij}^k\neq 0$, then $\epsilon_i\epsilon_j\epsilon_k=1$.\label{fusionlabel} 

\end{enumerate}

\end{proposition}

The proof is left as an exercise.  We remark that the sign $\epsilon_\alpha$ has appeared before under the name \textit{monodromy charge} \cite{SY,FRS}.
Using the signs $\epsilon_i$ of labels, we define a $\Z_2$-grading (on simple objects) as follows:
a simple object $X_i$ has a trivial grading or is in the \textit{local or trivial or even sector} $\mcB_0$ if $\epsilon_i=1$; Otherwise, it has a non-trivial grading or is in the \textit{twisted or defect or odd sector} $\mcB_1$.
Let ${I}_0$ be the subset of ${\Pi_\mcB}$ consisting of all labels in the trivial sector $\mcB_0$, and ${I}_1$ all labels in the defect sector $\mcB_1$.

\begin{proposition}\label{smatrixdecomp}

Let $(\mC,f)$ be a spin modular category, then

\begin{enumerate}

\item $\Pi_\mC={I}_0\coprod {I}_1$, and $f\in {I}_0$.

\item The tensor product respects the $\Z_2$-grading, $\mC = \mC_0\oplus \mC_1$. In particular, the action of $f$ on $\mcC$ by $\ot$ preserves the $\Z_2$-grading, and hence induces an action on $I_0$ and $I_1$.

\item If a simple object $\alpha$ is fixed by $f$ then $\alpha\in \mcC_1$ is a defect object. In particular, the action of $f$ restricted to ${I}_0$ is fixed-point free. 

\item If a simple object $\alpha$ is fixed by $f$, then for any $j\in {I}_1$, we have $s_{\alpha j}=0$.  If $s_{\alpha j}\neq 0$, then $j\in {I}_0$.

\item Let $I$ be a set of representatives of the orbits of the $f$-action on ${I}_0$, and $I_f=I_0\setminus I$. Partition the defect labels $I_1=I_{1n}\coprod I_{1f}$ into  non-fixed points and fixed points of the $f$-action, respectively. If the normalized $S$-matrix of $\mC$ is written in a $4\times 4$ block form indexed by $I, I_f, I_{1n}, I_{1f}$, then three of the $16$ blocks are $\bf{0}$, i.e., $s$ decomposes as follows:

$$s=\begin{pmatrix}
s_{II}&s_{II_f}&s_{II_{1n}}&s_{II_{1f}}\\
s_{I_fI}&s_{I_fI_f}&s_{I_fI_{1n}}& s_{I_fI_{1f}}\\
s_{I_{1n}I}&s_{I_{1n}I_f}&s_{I_{1n}I_{1n}}& \bf{0} \\
s_{I_{1f}I}& s_{I_{1f}I_f}&\bf{0} & \bf{0}
  \end{pmatrix}$$\label{blockform}

\item $f\cdot \alpha^*=(f\alpha)^*$, i.e., $f$ is compatible with duality or charge conjugation.

\end{enumerate}

\end{proposition}

\begin{proof}

(i):  Obvious from the definition.

(ii):  Obvious from (\ref{fusionlabel}) of Prop. \ref{fermionproperties}.

(iii):  By (\ref{thetasign}) of Prop. \ref{fermionproperties}, we have $\theta_\alpha=\theta_{f\alpha}=-\epsilon_\alpha\theta_\alpha$.  Hence $\epsilon_\alpha=-1$.

(iv):  By (\ref{smatrixsign}) of Prop. \ref{fermionproperties}, we have $s_{\alpha,j}=\epsilon_j s_{\alpha,j}$.  So if $j$ is in the defect sector, then $s_{\alpha,j}=0$.  But if $s_{\alpha,j}\neq 0$, then $\epsilon_j=1$, i.e., $j$ is in the trivial sector.

(v) and (vi): Obviously.

\end{proof}

\subsection{Fermionic Modular Categories}

Given a spin modular category $(\mC,f)$, then $\mcC=\mcC_0\oplus \mcC_1$, where $\mcC_i,i=0,1$ are the trivial and defect sectors, respectively.   Condensing $f$ results in a quotient category $\qcal$ of $\mC$.  The quotient category $\qcal$ encodes topological properties of the fermion system such as the ground state degeneracy of the system on the torus.  In the quotient, \lq\lq the fermion $f$ is condensed" because it is identified with the ground state represented by the tensor unit $\one$.  Naive fusion rules for $\qcal$ can be obtained by identifying objects in the orbits of the $f$-action as in Definition \ref{quotientfusionrule} below.  This idea goes back at least to M\"uger \cite{M3}: in his Proposition/Definition 2.15 where he obtains a tensor category from an idempotent completion of the category of $\Gamma$-modules in $\mcC$ for an algebra $\Gamma$.  In order to get a (linear) tensor category, he (implicitly) assumes $\Gamma$ is commutative, whereas our algebra object $\1\oplus f$ is not, so we do not obtain a fusion category (as in Example \ref{su2-6}.  It is an interesting question to formalize the quotient $\mcQ$ categorically, see \cite{Usher} for some progress.  This lack of a tensor structure, a braiding or a twist makes $\qcal$ unwieldy to work with.  Instead we will focus on some closely related categories: the two-fold \emph{covering theory} $(\mC,f)$ of $\qcal$, the trivial sector $\mcC_0\subset \mcC$ (for which $\mcC$ is also a 2-fold covering in a different sense) and the fermionic quotient $\qcal_0$ of $\mcC_0$.  The latter two reductions are motivated as follows: 1) in physical applications, sometimes we discard the defect sector $\mcC_1$ because the defect objects are not local with respect to the fermion $f$ and 2) the quotient $\qcal_0$ is better behaved than $\qcal$ (see Prop. \ref{smatrixproperties}).

\begin{definition}[see \cite{M3} Proposition/Definition 2.15]\label{quotientfusionrule}

Given a spin modular category $(\mcC,f)$ the object $\1\oplus f$ has a unique structure of an algebra (see the proof of Theorem 6.5 in \cite{KO}) : it is isomorphic to the (non-commutative) twisted group algebra $\mathbb{C}^w[\mathbb{Z}_2]$.  The following quotient $\qcal$ is called the \textit{fermionic  quotient} of $\mC$.   The objects of $\qcal$ are the same as $\mC$.  For two objects $x,y$ in $\qcal$, $\textrm{Hom}_{\qcal}(x,y)=\textrm{Hom}_{\mC}(x,y\otimes (\one \oplus f))$. Other structures such as braiding of $\mC$ will induce structures on $\qcal$.  The \textit{fermionic modular quotient} $\qcal_0$ of $\mcC_0$ is defined analogously.

\end{definition}

Let $[I_k], k=0,1$ be the orbit space of ${I}_k$ under the induced action of the fermion $f$.  Elements of $[I_k]$ are equivalence classes of labels in ${I}_k$, so the corresponding class of $i\in \Pi_\mcC$ will be denoted by $[i]$.  The label set of the fermionic quotient $\qcal$ is $[I_0]\cup [I_1]$, whereas the fermionic modular quotient $\qcal_0$ has label set $[I_0]$.  Given labels $[i],[j], [k] $, choose $i, j,k$ in $\Pi_\mcC$ covering $[i],[j],[k]$.  Then the naive fusion rules are $N_{[i][j]}^{[k]}=N_{ij}^k+N_{ij}^{fk}=\dim\Hom(X_i\ot X_j,X_k\ot(\1\oplus f))$.
  
Define ${D}_k^2=\sum_{i\in {I}_k} d_i^2,k=0,1$.
\begin{proposition}\label{smatrixproperties}
Let $(\CC,f)$ be a spin modular category.  Then:
\begin{enumerate}

\item ${D}_0^2={D}_1^2$.

\item $d_i=d_{fi}$.  Therefore, the quantum dimensions of labels descend to $[I_k],k=0,1$.

\item The braiding satisfies:

\begin{eqnarray}
 c_{j,fi}\cdot c_{fi,j}&=&\epsilon_j\Id_f\ot( c_{j,i}c_{i,j}),\\
 c_{fj,fi}c_{fi,fj}&=&\epsilon_i \epsilon_j(\Id_f\ot[(c_{i,f}^{-1}\ot \Id_j)(\Id_f\ot(c_{j,i}c_{i,j})(c_{i,f}\ot \id_j)]
\end{eqnarray} \label{braidingsign}

Therefore, pure braidings are well-defined on $[I_{0}]$, but ill-defined on $[{I}_1]$.  It follows that the $S$-matrix of $\mC$ descends to a well-defined matrix indexed by $[I_0]$, but does not descend to $[I_1]$.
\item The $T$-matrix of $\mC$ descends to a well-defined matrix indexed by $[I_1]$.  Although twists $\{\theta_i\}$ do not descend to $[I_0]$, double twists do descend to $[I_0]$.

\end{enumerate}

\end{proposition}

\begin{proof}

(i): By unitarity of the $S$-matrix,
$\sum_{j\in {I}}\tilde{s}_{0,j} \overline{\tilde{s}_{f,j}}=0$.  Since $\tilde{s}_{f,j}=\epsilon_jd_j$, $\sum_{j\in {I}}\tilde{s}_{0,j}\overline{\tilde{s}_{f,j}}=\sum_{j\in {I}}\epsilon_j d_j^2=0$.  The desired identity follows because
$\sum_{j\in {I}} \epsilon_j d_j^2={D}_0^2-{D}_1^2.$

(ii): We have $d_{fi}=d_fd_i=d_i$.

(iii): The first equation follows from funtoriality of the braiding and Proposition \ref{fermionproperties}(ii), and the second equation follows from the first.

(iv): Follows from (\ref{thetasign}) of Prop. \ref{fermionproperties}.

\end{proof}

By Prop. \ref{smatrixproperties}, we can define an $S$-matrix labeled by $[{I}_0]$.  To normalize correctly, we set $[s]_{[i],[j]}=2s_{i,j}$ for any $i,j\in {I}_0$.  Let $s=(s_{ij}), i,j\in [I_0]$, i.e., $s=2s_{II}$.  Notice that $s$ is symmetric since $s_{II}$ appears on the diagonal of the $S$-matrix of $\mcC$.  If we set $[D]^2=\sum_{i\in [I_0]} d_i^2$, then $[D]^2=\frac{1}{4} {D}^2=\frac{1}{2} {D}^2_0$.

\begin{theorem}

Given a spin modular category  $(\mC,f)$:

\begin{enumerate}

\item The matrix $[s]$ is unitary.

\item The Verlinde formula holds, i.e.,  $N_{[i][j]}^{[k]}=\sum_{r\in [I_0]}\frac{[s]_{i,r}[s]_{j,r}\overline{[s]_{k,r}}}{[s]_{0,r}}$ for any $[i],[j],[k]\in [I_0]$.

\end{enumerate}

\end{theorem}

\begin{proof}

(i): Given $i,j\in {I}_0$, we have $\sum_{k\in {I}} s_{ik}\overline{s_{kj}}=\delta_{ij}$, and $\sum_{k\in {I}} s_{if,k}\overline{s_{kj}}=
\sum_{k\in {I}} \epsilon_k s_{ik}\overline{s_{kj}}=\delta_{if,j}$.  If $j\neq i, fi$, then $\sum_{k\in {I}_0} s_{ik}\overline{s_{kj}}=\sum_{k\in {I}_1} s_{ik}\overline{s_{kj}}=0$.  Otherwise, we may assume $j=i\neq fi$.  Then $\sum_{k\in {I}_0} s_{ik}\overline{s_{ki}}=\sum_{k\in {I}_1} s_{ik}\overline{s_{k}}$ and $\sum_{k\in {I}_0} s_{ik}\overline{s_{ki}}=\frac{1}{2}$.  Since each $k\in [I_0]$ is covered by $2$ in ${I}_0$, we have $\sum_{k\in [{I}_0]} s_{ik}\overline{s_{kj}}=\frac{1}{4}$.  It follows that
$\sum_{k\in [{I}_0]} [s]_{ik}\overline{[s]_{kj}}=\delta_{ij}$.

(ii):  For any $i,j,k\in {I}_0$, we have
$$N_{ij}^k=
\sum_{r\in {I}_0}\frac{s_{i,r}s_{j,r}\overline{s_{k,r}}}{s_{0,r}}+\sum_{r\in {I}_1}\frac{s_{i,r}s_{j,r}\overline{s_{k,r}}}{s_{0,r}}.$$

Consider the same formulas for $N_{ij}^{fk}$.  For the first term
$\sum_{r\in {I}_0} \frac{s_{i,r} s_{j,r} \overline{s_{fk,r}} }{s_{0,r}}= \sum_{r\in {I}_0} $
$\epsilon_r \frac{s_{i,r} s_{j,r}\overline{s_{k,r}}} {\s_{0,r}}=\sum_{r\in {I}_0} \frac{s_{i,r}s_{j,r}\overline{s_{k,r}}} {s_{0,r}}.$ For the second term $\sum_{r\in {I}_1} \frac{s_{i,r} s_{j,r}\overline{s_{fk,r}}} {s_{0,r}}=
\sum_{r\in {I}_0}\epsilon_r \frac{s_{i,r} s_{j,r} \overline{s_{k,r}}} {s_{0,r}}=-
\sum_{r\in {I}_1}\frac{s_{i,r}s_{j,r}\overline{s_{k,r}}}{s_{0,r}}.$
Therefore, 
$$N_{ij}^k+N_{ij}^{fk}=2 \sum_{r\in {I}_0}\frac{s_{i,r}s_{j,r}\overline{s_{k,r}}}{s_{0,r}}=2(\sum_{r\in I_0}\frac{s_{i,r}s_{j,r}\overline{s_{k,r}}}{s_{0,r}}+\sum_{r\in I_0}\frac{s_{i,fr}s_{j,fr}\overline{s_{k,fr}}}{s_{0,r}})=4\sum_{r\in [I_0]}\frac{s_{i,r}s_{j,r}\overline{s_{k,r}}}{s_{0,r}},$$  which is the desired Verlinde formula.

\end{proof}

\subsection{Mapping class group representations}

A modular category gives rise to a unitary representation of the mapping class groups of the torus $T^2$, which is isomorphic to $SL(2,\Z)$.  A general quotient $\qcal$ of a spin modular category $\mcC$ is not a modular category, so we do not expect the existence of a representation of $SL(2,\Z)$. 

However, observe that $\qcal_0$, despite having no complete categorical description, has some of the data of a modular category: $\qcal_0$ has (naive) fusion rules and a unitary $S$-matrix obeying the Verlinde formula.
Moreover, the (normalized) $S$-matrix $s=2s_{II}$ and squared $T$-matrix $T^2$ are well-defined, and $s^4=I$.  A natural question is to ask if $s$ and $T^2$ combine to give a representation of a subgroup of $SL(2,\Z)$.

The subgroup $\Gamma_\theta$ of $SL(2,\Z)$ generated by $\mathfrak{s}=\begin{pmatrix}
0&-1\\
1 &0
  \end{pmatrix}$ and $\mathfrak{t}^2=\begin{pmatrix}
1&2\\
0 &1
  \end{pmatrix}$
 is isomorphic to the modular subgroup $\Gamma_0(2)$ consisting of matrices in $SL(2,\Z)$ that are upper triangular modulo $2$, via conjugation by $\begin{pmatrix}
 1 & 1\\0&2
 \end{pmatrix}$.  Projectively,  the images of $u$ and $v$ are independent so that as an abstract group $\Gamma_\theta/(\pm I)$ is generated by $\mathfrak{s},\mathfrak{t}^2$ satisfying $\mathfrak{s}^2=1$.  Therefore we have:
\begin{theorem}

The assignments $\mathfrak{s}\mapsto s$ and $\mathfrak{t}^2\mapsto T^2$ defines a projective representation of the group $\Gamma_\theta$ which does not come from a representation of $PSL(2,\Z)$ if the fermionic modular quotient $\qcal_0$ is not of rank=$1$.

\end{theorem}
We remark that ${s}^2$ and $T$ are well-defined on $[I_1]$, but since ${s}^2$ is the charge conjugation (permutation) matrix, this representation is not as interesting.

\subsection{Examples}

Spin modular categories that model fermionic quantum Hall states have well-defined fractional electric charges for anyons, i.e. another $\mathbb{Z}_n, n\geq 3$ grading beside the $\mathbb{Z}_2$ grading.  When a spin modular category ${\mC}$ comes from representations of an $N=2$ super conformal field theory, the sectors $\mcC_k,k=0,1$ are the Neveu-Schwartz (NS) and Ramond (R) sectors, respectively.

\begin{example}\label{mrcat}

The Moore-Read theory is the leading candidate for the fractional quantum Hall liquids at filling fraction $\nu=\frac{5}{2}$. The spin modular category of the Moore-Read theory is $Ising\times \Z_8$ with the fermion $f=\psi\otimes 4$.  The trivial NS sector consists of $\{\one\otimes i,\psi\otimes i\}$ for $i=$even and $\{\sigma \otimes i\}$ for $i=$odd.  Somewhat surprisingly, the rank=$6$ fermionic modular quotient theory can be given the structure of a linear monoidal category (see \cite[Appendix A.1.25]{Bon}) with labels $\{\one, \psi,\sigma, \bar{\sigma},\alpha,\bar{\alpha}\}$, where $\one, \psi$ are self-dual, $\sigma, \bar{\sigma}$ are dual to each other, and so are $\alpha, \bar{\alpha}$.  All fusion rules will follow from the following ones and obvious identities such as $\one x=x, xy=yx, \overline{xy}=\bar{y}\bar{x}$:

\begin{enumerate}

\item $\psi^2=1, \alpha\bar{\alpha}=1, \sigma\bar{\sigma}=1+\psi$

\item $\alpha^2=\psi, {\bar{\alpha}}^2=\psi, \sigma^2=\alpha+\bar{\alpha}, {\bar{\sigma}}^2=\alpha+\bar{\alpha}$

\item $\psi \sigma=\sigma, \psi \bar{\sigma}=\bar{\sigma}, \psi \alpha=\bar{\alpha}, \psi\bar{\alpha}=\alpha$

\item $\alpha \sigma=\bar{\sigma}, \alpha{\bar{\sigma}}=\sigma.$

If the labels are ordered as $\one, \sigma, \psi, \alpha, \bar{\sigma}, \bar{\alpha}$, then the $S$-matrix is

$$s=\frac{1}{2\sqrt{2}}\begin{pmatrix}
1 &\sqrt{2}& 1 & 1&\sqrt{2} & 1\\
\sqrt{2} & 0&-\sqrt{2} &i\sqrt{2} & 0 &-i\sqrt{2}\\
 1 &-\sqrt{2}& 1& 1&-\sqrt{2} & 1\\
 1 & i\sqrt{2} &1 &-1 &-i\sqrt{2} &-1\\
 \sqrt{2} & 0 &-\sqrt{2} &-i \sqrt{2} &0 &i\sqrt{2}\\
 1 &-i\sqrt{2} &1 &-1 &i\sqrt{2} &-1
  \end{pmatrix}$$

\end{enumerate}

\end{example}

Since this set of fusion rules comes from the subquotient $\qcal_0$ of a spin modular category, we expect there is a realization by a unitary fusion category without braidings.  Actually, the above fusion rules cannot be realized by any braided fusion category \cite{Bon}.

\begin{example}\label{su2-6}

Consider the spin modular category $SU(2)_6$.  The label set is ${I}=\{0,1,2,3,4,5,6\}$ and $6$ is the fermion.  Then ${I}_0=\{0,2,4,6\}$, ${I}_1=\{1,3,5\}$, $[I_0]=\{0,2\}$ and $[I_1]=\{1,3\}$.  $SU(2)_6$ is not graded for any $\mathbb{Z}_n, n\geq 3$.

Let $X_0=\one, X_2=x$ and define the fusion rules for the quotient as in (\ref{quotientfusionrule}), then we have:
$$ x^2=\one+2x.$$

It is known \cite{O2} that there are no fusion categories of rank=$2$ with fusion rules $x^2=\one+2x$, so $[I_0]$ cannot be the label set of a fusion category.  But there is a fermionic realization of the rank=$2$ category $\{1,x\}$ with $x^2=1+2x$ using solutions of pentagons with Grassmann numbers \cite{Cheng}.

The $S$-matrix as defined above is
$$s=\frac{1}{\sqrt{4+2\sqrt{2}}}
\begin{pmatrix}
1&1+\sqrt{2} \\
1+\sqrt{2}& -1
  \end{pmatrix}$$

Note although Verlinde formulas do give rise to the above fusion rules, this unitary matrix is not the modular $s$-matrix of any rank=$2$ modular category.
\end{example}

\begin{example}
Laughlin fractional quantum Hall states at filling fraction $\nu=\frac{1}{Q}$, $Q$=odd, has $Q$ different anyons labeled by $r=0,1,\ldots, Q-1$.  Note $Q=1$ is an integer quantum Hall state.  The conformal weight of anyon $r$ is $h_r=\frac{r^2}{2Q}$.

The covering spin modular category is the abelian UMC $Z_{4Q}$ labeled by $a=0,1,\cdots, 4Q-1$.  The twist of the object $a$ is $\theta_a=\frac{a^2}{8Q}$.  Its charge is $q_a=\frac{a}{2Q}$.  The fermion corresponds to $f=2Q$.  The double braiding of $i,j$ is $\lambda_{ij}=e^{\frac{2\pi \ii ij}{4Q}}=\tilde{s}_{ij}$.  It follows that $\tilde{s}_{i,f}=(-1)^i$, hence ${I}_0$ consists of all even labels, while ${I}_1$ are the odd labels.
\end{example}

We end this section with:

\begin{question}

Given a spin modular category $(\mC,f)$, are the following true?

\begin{enumerate}

\item The fermion $f$ has no fixed points if and only if $N_{ij}^k\cdot N_{ij}^{fk}=0$ for all $i,j,k$.

\item If $f$ has no fixed points, then  $(\mC,f)$ has a $\mathbb{Z}_n, n\geq 3$ grading.

\end{enumerate}

\end{question}

\section{Super-modular categories}

Let $\mcB$ be a braided fusion category, and $\mathcal{D} \subset \mcB$ a fusion subcategory. The M\"{u}ger
centralizer $C_\mcB(\mathcal{D})$ of $\mathcal{D}$ in $\mcB$ is the fusion subcategory generated by $Y\in Ob(\mcB)$ such that $c_{Y,X}\circ c_{X,Y}=\id_{X \otimes Y}$ for any $X$ in  $\mathcal{D}$. The M\"{u}ger center of $\mcB$ is the symmetric fusion subcategory $\mathcal{Z}_2(\mcB):=C_\mcB(\mathcal{B})$. The objects of $\mathcal{Z}_2(\mcB)$ are called \emph{transparent}, and we sometimes use the shorthand $\mcD^\prime$ for $C_\mcB(\mcD)$ when no confusion can arise.


\begin{definition}\label{super-modular}
A URFC $\mcB$ is called \textbf{super-modular} if its M\"uger center $\mathcal{Z}_2(\mcB)\cong \sVec$, i.e. every non-trivial transparent simple object is isomorphic to the same fermion.
\end{definition}
Without the unitarity and sphericity assumptions, braided fusion categores with M\"uger center $\sVec$ (as a symmetric fusion category) are called slightly degenerate modular categories in \cite{DGNO1}.

The trivial sector $\mcC_0$ of a spin modular category $(\mcC,f)$ is a super-modular category.  It is not known if all super-modular categories arise this way and we conjecture that it is indeed so and provide evidence in this section.  Most of the results in the previous section proved for the trivial sector $\mcC_0$
of a spin modular category $(\mcC,f)$ can be proved directly for super-modular categories.

If $\CC$ is a UMC, then $\sVec\boxtimes\CC$ is super-modular.  If $\mcB\cong \sVec\boxtimes \CC$ with $\CC$ modular, we will say $\mcB$ is \emph{split super-modular}, and otherwise \emph{non-split super-modular}. Observe that a super-modular category is split if, and only if, it is $\Z_2$-graded with the corresponding trivial component modular.  In particular $\sVec$ is a split super-modular category since $\Vec$ is modular.

\begin{theorem}\label{theor on split super-mod}
Let $(\CC,f)$ be a spin modular category and $\mcC_0$ be the associated super-modular subcategory. Then the following are equivalent 

\begin{itemize}
\item[(i)] $\mcC_0$ is split super-modular.
\item[(ii)] $\mathcal{C}_0$ contains a modular subcategory of dimension $\dim(\mathcal{C}_0)/2.$
\item[(iii)] $\mathcal{C}$ contains a modular subcategory of dimension four that contains $f$.
\end{itemize}
\end{theorem}
\begin{proof}
Obviously $(i)$ implies $(ii).$

Assume $(ii)$. Let $\mathcal{D}\subset \mcC_0$ a modular category with $\dim(\mcC_0)=2\dim(\mathcal{D})$.
Since $\mathcal{D}\subset \mathcal{C}$ and $\CC$ is modular, it follows from \cite[Theorem 4.2]{M2} that $\CC=\mathcal{D}\boxtimes C_{\CC}(\mathcal{D})$, where $C_{\CC}(\mathcal{D})$ is modular and 
$$
    \dim(C_{\CC}(\mathcal{D}))=\frac{\dim(\CC)}{\dim(\mathcal{D})}\\
    = \frac{2\dim(\mcC_0)}{\dim(\mcC_0)/2}=4.
$$
Since $\mathcal{D}\subset \mcC_0$, we have that $\langle f\rangle = C_{\CC}(\mcC_0)\subset C_{\CC}(\mathcal{D})$. Hence $(ii)$ implies $(iii)$.  

Assume $(iii)$. Let  $\mathcal{A}\subset \CC$ be a modular subcategory with  $f\in \mathcal{A}$. Then $\CC=C_{\CC}(\mathcal{A})\boxtimes \mathcal{A}$ and $C_{\CC}(\mathcal{A})\subset \mcC_0$.  Since $C_{\CC}(\mathcal{A})\boxtimes \langle f\rangle \subset \mcC_0$ and  
\begin{equation*}
  \dim(C_{\mathcal{C}}(\mathcal{A})\boxtimes\langle f\rangle)=2\dim(C_{\mathcal{C}}(\mathcal{A}))=\dim(\mcC_0)
\end{equation*}
we have that $\mcC_0=C_{\CC}(\mathcal{A})\boxtimes \langle f\rangle$. Hence $\mcC_0$ is split super-modular.
\end{proof}

Let $G$ be a finite group and $w \in Z^3(G,\mathbb{C}^*)$. Recall the definition of $\beta_x(y,z)$ given in equation \eqref{gamma-defn}.



\begin{definition}(\cite{MR2552301})\label{inv-bichar}
Let $H,K$ be normal subgroups of $G$ that centralize each other. 
An $w$-{\em bicharacter} is a function $B:K\times H\rightarrow \mathbb{C}^{\times}$ such that
\begin{eqnarray*}
  \mbox{(i) }B(x,yz) &=& \beta_x^{-1}(y,z) B(x,y)B(x,z) \ \mbox{ and}\\
  \mbox{(ii) }B(sx,y) &=& \beta_y(s,x) B(s,y)B(x,y)
\end{eqnarray*}
for all $s,x\in K$, $y,z\in H$.

An $w$-bicharacter $B$ is called $G$-{\em invariant} if 
$$
   B(x^{-1}kx,h) = \frac{\beta_k(x,h)\beta_k(xh,x^{-1})}{\beta_k(x,x^{-1})}
                       B(k,xhx^{-1})
$$
for all $x,y\in G$, $h\in H$, $k\in K$.
\end{definition}

We recall the classification of fusion subcategories of $\operatorname{Rep}(D^w(G))$ given in \cite[Theorem 1.2]{MR2552301}. The fusion subcategories of  $\operatorname{Rep}(D^w(G))$  are in bijection with triples $(K, H, B)$ where $K$, $H$ are normal subgroups of $G$ centralizing
each other and $B : K \times H \to \mathbb{C}^*$ is a $G$-invariant $w$-bicharacter. The fusion subcategory associated a triple $(K, H, B)$ will be denoted $\mathcal{S}(K, H, B)$. 

\begin{remark}\label{rmk paper Deepak}
The following are some results from \textit{loc. cit.} that we will need.

\begin{itemize}
\item The dimension of $\mathcal{S}(K, H, B)$ is $|K|[G:H]$ (see \cite[Lemma 5.9]{MR2552301}).
\item $\mathcal{S}(K, H, B)\subset \mathcal{S}(K', H', B')$ if and only if $K\subset K', H'\subset H$ and
$B|_{K\times H'} = B'
|_{K\times H'}$, (see \cite[Proposition 6.1]{MR2552301}).
\item $\mathcal{S}(K, H, B)$ is modular 
if and only if $HK = G$ and the symmetric bicharacter $BB^{\operatorname{op}}|_{(K \cap H)\times(K\cap H)}$
is nondegenerate (see \cite[Proposition 6.7]{MR2552301}).
\end{itemize} 
\end{remark}

Recall that by Proposition \ref{fermions in D(G)} fermions in $\operatorname{Rep}(D^w(G))$ are in correspondence with pairs $(\eta,z)$, where $z$ is central element of order two and $\eta:G\to \mathbb{C}^*$ is a map satisfying some conditions, see \textit{loc. cit.} Applying Theorem \ref{theor on split super-mod}, the following proposition provides necessary and sufficient group-theoretical conditions in order that a super-modular category obtained from a spin modular twisted Drinfeld double be non-split. 

\begin{proposition}
Let $f$ be a fermion in $\operatorname{Rep}(D^w(G))$ with associated data $(\eta,z)$. The modular subcategories of $\operatorname{Rep}(D^w(G))$ of dimension $4$ containing $f$ correspond to: 
\begin{itemize}
    \item Subgroups $H\subset G$ such that $G=H\times \langle z\rangle$. The modular category associated to $H$ is  $\mathcal{S}(\langle z\rangle, H, B_\eta)$, where $B_\eta(z,x)=\eta(x)$ for all $x\in H$.
 \item    Pairs $(K,B)$, where $K\subset G$ is a central subgroup of order four containing $z$ and $B:K\times G\to \mathbb{C}^*$ is a $G$-invariant $w$-bicharacter such that 
\begin{enumerate}
\item $\eta(x)=B(z,x),$ for all $x\in G$.
\item The symmetric bicharacter $BB^{\operatorname{op}}:K\times K\to U(1)$ is nondegenerate.
\end{enumerate}
The modular category associated to $(H,B)$ is  $\mathcal{S}(H, G, B)$.
\end{itemize}
\end{proposition}
\begin{proof}

Let $f\in \operatorname{Rep}(D^w(G))$ be a fermion with associated data $(\eta,z)$, see Proposition \ref{fermions in D(G)}. The fusion subcategory generated by $f$ corresponds to $\langle f\rangle =\mathcal{S}(\langle z \rangle, G, B_\eta)$, where $B_\eta(z,x)=\eta(x)$ for all $x\in G$.

Let $\mathcal{S}(K, H, B)$ be a modular subcategory of $\operatorname{Rep}(D^w(G))$ of dimension $4$ containing $f$. Using the results cited in Remark \ref{rmk paper Deepak} we have
\begin{itemize}
\item[(a)] $|K|[G:H]= 4$
\item[(b)] $\langle f\rangle \subset K$
\item[(c)] $KH=G$.
\end{itemize}
The conditions $(a)$ and $(b)$ imply that there are only two possibilities:

\begin{itemize}
\item[(i)] $K=\langle z\rangle$ and $[G:H]=2$. 
\item[(ii)] $H=G$ and $K$ is a central subgroup of order four.
\end{itemize}
In the case that $K=\langle z\rangle$ and $[G:H]=2$. Condition $(c)$ implies that if $z\notin H$, then $G\cong H\times K$.
\end{proof}
\subsection{Braided fusion categories with transparent fermions}

The following is a structure theorem for unitary ribbon fusion categories $\mcB$. Transparent objects of  $\mcB$ form a symmetric fusion subcategory ${\mcB^\prime}$.  By a theorem of Deligne, every symmetric fusion category is equivalent to the representation category of a pair $(G,z)$, where $G$ is a finite group and $z$ is a central element of $G$ of order $\leq 2$ (see \cite{O1}).  
 $\mcB^\prime$ is Tannakian if and only if $z=1$.  Recall that a ribbon fusion category with a non-Tannakian M\"uger center need not have a transparent fermion (Proposition \ref{prop: symmetric contains fermion}).

\begin{lemma}\cite[Theorem 2]{Sawin}
There is a Tannakian subcategory $\mcS\cong\Rep(G)$ of a URFC $\mcB$ such that the de-equivariantization $\mcB_G$ is either modular or super-modular $\mcB_T$.
\end{lemma}

It follows that if $\mcB_G$ is not modular, there is an exact sequence:
$ 1\rightarrow \operatorname{Rep}(G) \rightarrow \mcB \rightarrow \mcB_T\rightarrow 1,$
where $\mcB_T$ is super-modular.  So a URFC is a twisted product of a Tannakian category and a super-modular category, therefore, a \lq\lq braided" equivariantization of a super-modular category.

\begin{proposition}\label{example}
Let $(\mathcal{B},f)$ be a super-modular category and $*:\underline{G}\to \operatorname{Aut}_{\otimes}^{br}(\mathcal{B})$ an action by a finite group $G$ such that the restriction of the $G$-action to $\langle f \rangle$ is trivial. Then the equivariantization $\mathcal{B}^G$ is pre-modular category with $\mathcal{Z}_2(\mathcal{B}^G)=\operatorname{Rep}(G)\boxtimes \sVec$. Moreover, every pre-modular category with a transparent fermion is constructed in this way. 
\end{proposition}
\begin{proof}
That $\mcB^G$ has the desired properties follows from definition of equivariantization. To prove every pre-modular category with a transparent fermion is constructed in this way, let $\mathcal{D}$ be a braided fusion category and $f\in \mathcal{Z}_2(\mathcal{D})$ a transparent fermion. By Proposition \ref{prop: symmetric contains fermion}, $\mathcal{Z}_2(\mathcal{D})=\operatorname{Rep}(G)\boxtimes \sVec$ as braided fusion categories. Then the algebra $\mathcal{O}(G)$ of functions on $G$ is a commutative 
algebra in $\mathcal{Z}_2(\mathcal{D})\subset \mathcal{D}$. The category $\mathcal{D}_G$ of
left $\mathcal{O}(G)$-modules in $\mathcal{D}$ is a braided fusion category, called de-equivariantization of $\mathcal{D}$ by $\operatorname{Rep}(G)$, see \cite{DGNO1} for more details. Moreover, the
free module functor $\mathcal{D} \to  \mathcal{D}_G,  Y\mapsto \mathcal{O}(G)\otimes Y$ is  a surjective braided functor. Hence   $\mathcal{D}_G$ is a super-modular category with fermion object $\mathcal{O}(G)\otimes f$.

By \cite[Theorem 4.4]{DGNO1}, equivariantization and de-equivariantization are mutually inverse processes. The group $G$ acts on $\mathcal{O}(G)$  (by right translations) viewed as an
algebra in $\Rep(G)$. Then $G$ acts on the super-modular category $\mathcal{D}_G$. In particular the action of $G$ on the transparent fermion $\mathcal{O}(G)\otimes f\in \mathcal{D}_G$ is trivial.
\end{proof}

\begin{example}
Let $(\mathcal{B},f)$ be a non-split super-modular category and $F:\mathcal{B}\to \mathcal{B}$ a non-trivial braided autoequivalence such that $F(f)\cong f$. If $F$ has order $n$, it defines a non trivial group homomorphism from $\mathbb{Z}_n$  to the group of braided autormorphisms of $\mathcal{B}$. This group homomorphisms lifts to a categorical action of $\mathbb{Z}_n$ on $\mathcal{B}$ if and only if a certain third cohomology class $O(F)\in H^3(\mathbb{Z}_n, \operatorname{Aut}_\otimes (\operatorname{Id}_\mathcal{B}))$ is zero, see \cite[Theorem 5.5 and Corollary 5.6]{Gal}.
Since $H^3(\mathbb{Z}, \operatorname{Aut}_\otimes (\operatorname{Id}_\mathcal{B}))=0$, even if $O(F)\neq 0$, there is a group epimorphism $p:\mathbb{Z}_m\to \mathbb{Z}_n$ such that $p^*(O(F))=0$. Thus the group $\mathbb{Z}_m$ acts non-trivially on $\mathcal{B}$. By 
Proposition \ref{example} the equivariantization 
$\mathcal{B}^{\mathbb{Z}_m}$ is
a premodular category with $\mathcal{Z}_2(\mathcal{B}^{\mathbb{Z}_m})=\operatorname{Rep}(\mathbb{Z}_m)\boxtimes\sVec$.  
\end{example}

The $S$ and $T$ matrices of a super-modular category have the following special form:
\begin{theorem}\label{sproduct}
 If $\mcB$ is super-modular, then $\tilde{s}_\mcB=\hat{S}\otimes \tilde{s}_{\sVec}$ and $T_\mcB=\hat{T}\otimes T_{\sVec}$ for some invertible matrices $\hat{S}$ and $\hat{T}$.
\end{theorem}
\begin{proof}
 Suppose $\mcB$ is super-modular with fermion $f$.  Since $\mcB^\prime=\sVec$, we have $\tilde{s}_{X,f}=d_X$ for all simple $X$.  Moreover, there is a (non-canonical) partition of the simple objects into two sets: $X_0=\mathbf{1},X_1,\ldots,X_r,f\otimes X_0=f,f\otimes X_1,\ldots,f\otimes X_r$, since $X\otimes f\not\cong X$ for any $X$.    The balancing equation gives us:

$$-\theta_Xd_X=\tilde{s}_{X,f}\theta_X\theta_f=d_{f\otimes X}\theta_{f\otimes X}=d_X\theta_{f\otimes X}.$$
Thus $\theta_X=-\theta_{f\otimes X}$, and $T_\CC=\hat{T}\otimes T_{\sVec}$.
Now we just need to show that $\tilde{s}_{X,X}=\tilde{s}_{f\otimes X,f\otimes X}=\tilde{s}_{X,f\otimes X}$ for all simple objects $X$ so that $\tilde{s}_\mcB=\hat{S}\otimes \tilde{s}_{\sVec}$.  
Fix $X$, and suppose that $X\otimes X^{*}=\sum_{Y\in A}Y$ for some (multi-)set $A$.  This implies that $f\otimes X\otimes X^*=\sum_{Y\in A} f\otimes Y$, and $f\otimes Y$ is simple.
Now $\tilde{s}_{X,X}=\frac{1}{\theta_X^2}\sum_{Y\in A} d_Y\theta_Y$.  Computing:
$$\tilde{s}_{X,f\otimes X}=-\frac{1}{\theta_X^2}\sum_{Y\in A} d_{f\otimes Y}\theta_{f\otimes Y}=-\frac{1}{\theta_X^2}\sum_{Y\in A} d_Y(-\theta_Y)=\tilde{s}_{X,X}.$$  Since $f$ is transparent, we also have $\tilde{s}_{f\otimes X,f\otimes X}=\tilde{s}_{X,X}$.

\end{proof}
The following is an immediate consequence:
\begin{corollary}\label{gauss sum is zero}
If $\mcB$ is a super modular category   \[\tau^\pm(\mcB):=\sum_{X\in \operatorname{Irr}(\mcB)} \theta_X^{\pm}d_X^2=0.\]
\end{corollary}

\subsection{Super-modular categories from quantum groups} 

Quantum groups at roots of unity yield unitary modular categories via ``purification" of representation categories (see \cite[Section XI.6]{TuraevBook} and \cite{Rsurv}).  By taking subcategories we obtain several non-split super-modular categories.

The modular category $SU(2)_{4m+2}$ obtained as a semisimple subquotient of the category of representations of the quantum group $U_q\mathfrak{sl}_2$ at $q=e^{\pi i/(4m+4)}$ has rank $4m+3$, with simple objects 
labeled $X_0=\mathbf{1},X_1,\ldots,X_{4m+2}$, (cf. \cite[Example 3.3.22]{BK}).  The $S$- and $T$-matrices are given by:

$\tilde{s}_{i,j}=\frac{\sin((i+1)(j+1)\pi/(4m+4))}{\sin(\pi/(4m+4))}$ and $t_{j,j}=e^{\pi i(j^2+2j)/(8m+8)}$. The object $X_{4m+2}$ is the only non-trivial invertible object and hence the universal grading group of $SU(2)_{4m+2}$ is $\Z_2$.

\begin{lemma}
 The subcategory,  $PSU(2)_{4m+2}$, of $SU(2)_{4m+2}$ generated by the $2m+2$ simple objects with even labels: $X_0=\mathbf{1},X_2,\ldots,X_{4m+2}$ is non-split super-modular.
\end{lemma}
\begin{proof}  We must show that the M\"uger center of $PSU(2)_{4m+2}$ is isomorphic to $\sVec$.  Since the M\"uger center is always a symmetric (and hence integral) category we first observe that the only non-trivial object with integral dimension is $X_{4m+2}$, in fact $\dim(X_{4m+2})=1$. It is routine to check that $\tilde{s}_{4m+2,2j}=\dim(X_{2j})$ and that $\theta_{4m+2}=e^{\pi i(4m+2)(4m+4)/(8m+8)}=-1$.  To see that $PSU(2)_{4m+2}$ is non-split super-modular observe that if $\CC$ were a modular subcategory of $PSU(2)_{4m+2}$ with rank $m+1$ then $SU(2)_{4m+2}$ would factor as a Deligne product of two modular categories.  But $m+1$ does not divide $4m+3$, so this is impossible.
\end{proof}

Observe that for $m=0$ we recover $\sVec=PSU(2)_2$.

A $2$-parameter family of non-split super-modular categories can be obtained as subcategories of $SO(N)_r$ for $N,r$ both odd, i.e. the modular category obtained from $U_q\mathfrak{so}_N$ with $q=e^{\frac{\pi i}{2(r+N-2)}}$.  Let $PSO(N)_r$ be the subcategory with simple objects labeled by the highest weights of $SO(N)_r$ with integer entries.   Identifying $SU(2)_{4m+2}$ with $SO(3)_{2m+1}$ the examples above can be made to fit into this larger family.  Setting $N=2s+1$ and $r=2m+1$ we compute the rank of $SO(2s+1)_{2m+1}$ to be $\frac{3s+4m}{s+m}\binom{s+m}{s}$, while the rank of $PSO(2s+1)_{2m+1}$ is $2\binom{s+m}{s}$ (here one uses the combinatorial methods described in \cite{Rsurv}).  The object $f$ in $PSO(2s+1)_{2m+1}$ labelled by the weight vector $r\Lambda_1=(r,0,\ldots,0)$ is a fermion, and $\otimes$-generates the M\"uger center of $PSO(2s+1)_{2m+1}$, which can be explicitly shown as in the $PSU(2)_{4m+2}$ case.  To see that $PSO(2s+1)_{2m+1}$ cannot be split supermodular observe that 1/2 the rank of $PSO(2s+1)_{2m+1}$ does not divide the rank of $SO(2s+1)_{2m+1}$, so $PSO(2s+1)_{2m+1}$ cannot factor as $\sVec\boxtimes \CC$ for some modular category $\CC$.

\subsection{The Modular Closure Conjecture}

\begin{definition}

\begin{enumerate}
    \item Let $\mcB$ be a URFC.  A modular category $\CC \supset \mcB$ is called a \textit{minimal modular extension or modular closure} of $\mcB$ if $D_{\mathcal{C}}^{2}=D_{\mcB^\prime}^{2}D_{\mcB}^{2}$.
    \item Two modular extensions $\CC_1 \supset \mcB$ and $\CC_2 \supset \mcB$ are equivalent if there is a  braided  equivalence $F: \CC_1 \to \CC_2$ such that $F|_{\mcB}=\operatorname{Id}_{\mcB}$. 
\end{enumerate}
\end{definition}
A minimal modular extension of a super-modular category $\mcB$ is a spin modular category $(\mC, f)$ with the fermion $f$ being the transparent one in $\mcB$. 
\subsubsection{Counterexamples to the modular closure conjecture}

Recall from \cite{M2}
\begin{conjecture}
  Let $\mcB$ be a URFC category, then there exists a UMC ${\mcC}$ and a full and faithful tensor functor $I:\mcB\to{\mcC}$ such that $D_{\mcC}^{2}=D_{\mcB}^{2}D_{\mcB^\prime}^{2}$.
 \end{conjecture}
 
M\"uger's modular closure conjecture as above in full generality does not hold.  Unpublished counterexamples due to Drinfeld exist \cite{OE}.  A general method for constructing counterexamples is the following:

Let $G$ be a finite group acting by braided-automorphisms on a modular category $\mcB$, $\rho:BG\to B\operatorname{Aut}_{br}(\mcB)$. Then $\mcB^{G}$ is again braided and its M\"{u}ger center is $\mathrm{Rep}(G)$. Now suppose that there exists a minimal modular extension $\mcB^{G}\subset \mcM$, then the de-equivariantization $\mcM_{G}$ is a faithful $G$-crossed modular category that corresponds to a map $BG\to BPic(\mcB)$ and it is a lifting of the $G$-action on $\mcB$. In other words, $\mcB^{G}$ admits a minimal modular extension if and only if $\rho$ admits a gauging. One can compute the obstruction explicitly in some cases. For instance, if $\mcB=\Vec_{A}$, and the modular structure is given by a bicharacter, then the obstruction is the cup product \cite{CGPW,Getal}.

Drinfeld proved that the obstructions in the following cases are nonzero:
\begin{itemize}
    \item  $G=(\Z_{2}\times\Z_{2})$, $\mcB=\mathrm{Sem}$, $\alpha\in H^{2}(G,\Z_2)$ corresponds to the Heisenberg group.
    \item $G=\Z_{p}\times \Z_{p}$, $\mcB=\Vec_{\Z_{p}}$ with the canonical modular structure $\alpha \in H^{2}(\Z_{2}^{2},\Z_{2})$ corresponding to an extensions non-isomorphic to the Heisenberg group.
\end{itemize}

\subsection{The 16-fold Way Conjecture}

A super-modular category models the states in the local sector of a fermionic topological phase of matter.  In physics, gauging the fermion parity should result in modular closures of super-modular categories by adding the twisted sectors.  In two spatial dimensions, gauging the fermion parity seems to be un-obstructed.  

\begin{conjecture}\label{16foldway}
Let $\mcB$ be super-modular.  Then $\mcB$ has precisely 16 minimal unitary modular extensions up to ribbon equiavence.
\end{conjecture}
In fact, in \cite{LKW} it is shown that if $\mcB$ has one minimal modular extensions then it has precisely 16.

\begin{lemma}
Suppose $\mcB$ is super-modular, and $\CC$ is a minimal modular extension of $\mcB$.  Then $\CC$ is faithfully $\Z_2$-graded with $\CC_0=\mcB$.
\end{lemma}

\begin{proof}  Since $\sVec\subset\CC$ and $\CC$ is modular, $\CC$ is faithfully $\Z_2$-graded, with trivial component $\CC_0=\sVec^\prime$.  Since $\mcB^\prime=\sVec$, we have $\mcB\subset\CC_0$.  Since $\dim(\mcB)=\dim(\CC_0)$, the proof is complete.
\end{proof}

The following result due to Kitaev \cite{kitaev} is the 16-fold way for free fermions:
\begin{prop} $\sVec$ has precisely $16$ inequivalent minimal unitary modular closures $SO(N)_1$ for $1\leq N\leq 16$, where $N=1$ denotes the Ising theory and $N=2$ the $U(1)_4$-cyclic modular category.  They are distinguished by their multiplicative central charges, which are
$e^{2\pi i \nu/16}$ for $1\leq \nu=N\leq 16$.
\end{prop}

In what follows we will denote $SO(N)_1$ by $\mS_\nu$ with $\nu=N$.
Kitaev's result immediately implies Conjecture \ref{16foldway} holds for split super-modular categories:

\begin{corollary}
If $\CC$ is modular then $\CC\boxtimes\sVec$ has precisely 16 inequivalent minimal modular closures.
\end{corollary}
\begin{proof}
  Clearly if $\mS_\nu$ is a minimal modular closure of $\sVec$ then $\CC\boxtimes\mS_\nu$ is a minimal modular closure of $\CC\boxtimes\sVec$.  On the other hand, if  $\mcd$ is a minimal modular closure of $\CC\boxtimes\sVec$ then $\mcd\cong\CC\boxtimes C_\mcd(\CC)$ with $C_\mcd(\CC)$ by \cite[Theorem 4.2]{M2}.  Thus $C_\mcd(\CC)$ is a minimal modular closure of $\sVec$ and hence $\mcd\cong\CC\boxtimes\mS_\nu$ for some $\mS_\nu$.
\end{proof}

\subsection{Witt class 16-fold way}

Witt equivalence for modular categories and the Witt group $W$ are defined in \cite[Section 5.1]{DMNO}. Super-Witt equivalence and the super-Witt group $sW$ are defined in \cite[Section 5.1]{DNO}.  The following two Theorems \ref{16ways} and \ref{allwitt} imply that if a super-modular category has one minimal modular extension then it has $16$ up to Witt equivalence (cf. \cite[Theorem 5.3]{LKW}).

\begin{theorem}\label{16ways}  
Let $\mcB$ be a super-modular category with a minimal modular extension $\mcC$ and transparent fermion $f$. Furthermore, let $e$ be a generator for $\sVec$ and $\mS_{\nu}$ and $\mS_{\mu}$ two inequivalent minimal modular extensions of $\sVec$. Then
\begin{itemize}
    \item[(i)] $(f, e)\in \mcC\boxtimes \mS_{\nu}$ generates a Tannakian subcategory, $\mcE\cong\Rep(\Z_{2})$.
    \item[(ii)] $\mcC_{\nu}:=[(\mcC\boxtimes\mS_{\nu})_{\mcE}]_0$ is a minimal modular extension of $\mcB$, with multiplicative central charge the same as that of $\mcC\boxtimes\mS_{\nu}$.
    \item[(iii)] $\mcC_{\mu}$ and $\mcC_{\nu}$ are Witt inequivalent, and hence inequivalent.
\end{itemize}
\end{theorem}

\begin{proof}
  It follows immediately from the definition of the Deligne product that $\mcC\boxtimes\mS_{\nu}$ is modular, and that$(f, e)$ generates a Tannakian subcategory, $\mcE\cong\Rep(\Z_2)$.
  In particular, $(\CC\boxtimes \mS_\nu)_{\EE}$ is $\Z_2$-crossed braided with modular trivial component $\mcC_\nu$ by \cite[Proposition 4.56(i)]{DGNO1}. Applying \cite[Proposition 4.26 and Corollary 4.28]{DGNO1} we find that $\dim(\mcC_\nu)=\dim(\mcC)$. By \cite{CGPW}, the multiplicative central charge can be computed as: $\xi(\mcC_\nu)=\xi(\mcC\boxtimes\mS_\nu)$, which is $\xi(\mcC)e^{\pi i\nu/8}$.
  So to prove (ii), it remains to show that $\mcB$ is a ribbon subcategory of $\mcC_{\nu}$. By \cite[Proposition 4.56(ii)]{DGNO1}, $\mcC_{\nu} = (\mcE')_{\Z_2}$, while the definition of $\mcE$ gives 
  \begin{align*}
    \mcE' = (\mcB \boxtimes \sVec) \oplus (\mcC_1 \boxtimes (\mcC_\nu)_1),
  \end{align*}
  where $\mcC_1$ and $(\mcC_\nu)_1$ are the odd gradings of $\mcC$ and $\mcC_\nu$ respectively. Since $(\mcB \boxtimes \sVec)_{\mathbb Z_2} = \mcB$ and the de-equivariantization respects the grading, (ii) follows.
  
  Finally, suppose $\mS_{\mu}$ is a minimal modular extension of $\sVec$ that is inequivalent to $\mS_{\nu}$. Then $\mS_{\mu}$ and $\mS_{\nu}$ have distinct (multiplicative) central charges. So, by \cite[Remark 6.17]{DGNO1}, it follows that $\mcC_{\nu}$ and $\mcC_{\mu}$ have inequivalent central charges. Thus (iii) follows from \cite[Lemma 5.27]{DMNO}
\end{proof}

\begin{theorem}\label{allwitt} If $\mcB$ is super-modular, then every minimal modular closure of $\mcB$ is Witt-equivalent to one of the extensions obtained in Theorem~\ref{16ways}.
\end{theorem}
\begin{proof}
  Let $\mcC$ be a minimal modular closure of $\mcB$. Then Witt class $[\mcC]_W$ is sent to the super-Witt class $[\mcC \boxtimes \sVec]_{sW}$ under the canonical homomorphism $g:W \to sW$ defined in \cite[Section 5.3]{DNO}. By \cite[Proposition 5.14]{DNO} we know that the the kernel of $g$ consists of the Witt classes represented by modular closures of $\sVec$. So by Theorem~\ref{16ways}(iii), it suffices to show that $\mcC\boxtimes \sVec$ and $\mcB$ are super-Witt equivalent. To this end, let $\mcE\subset\mcC\boxtimes\sVec$ be the Tannakian category described in the previous theorem. By \cite[Proposition 5.3]{DNO},
  \begin{align*}
    [\mcC \boxtimes \sVec]_{sW} = [((\mcC \boxtimes \sVec)_{\Z_2})_0]_{sW}.
  \end{align*}

  Finally, by \cite[Proposition 4.56(ii)]{DGNO1}, we have
  \begin{align*}
    ((\mcC \boxtimes \sVec)_{\Z_2})_0 = (\mcE')_{\Z_2} = (\mcB \boxtimes sVec)_{\Z_2} = \mcB.
  \end{align*}
\end{proof}

\subsection{Zested extensions of a super-modular category}
\subsubsection{\texorpdfstring{$G$}{G}-grading of modular categories}

It was proved in \cite[Theorem 3.5]{GN} that any fusion category $\mcC$ is naturally graded by a group $U(\mcC)$, called the universal grading group of $\mcC$, and the adjoint subcategory $\mcC_{ad}$ (generated by all subobjects of $X^*\ot X$, for all $X$) is the trivial component of this grading. Moreover, any other faithful grading of $\mcC$ arises from a quotient of $U(\mcC)$ \cite[Corollary 3.7]{GN}.

For any abelian group $A$, let denote $\widehat{A}$ the abelian group of linear complex characters.
For a braided fusion category, there is group homomorphism $\phi:U(\mcC) \to \widehat{G(\mcC)}$, roughly defined as follows: For $g\in G(\mcC)$ and $i\in\operatorname{Irr}(\mcC)$ the double braiding $c_{i,g}c_{g,i}$ is an isomorphism on the simple object $g\otimes i$, and hence a scalar map $\phi(i,g)\id_{g\otimes i}$ .  It can 
be shown that for each $i$, $\phi(i,-)$ is a character (and is related to the monodromy charge of \cite{SY,FRS}).  Therefore we obtain a multiplicative map $\phi:K_0(\mcC) \to \widehat{G(\mcC)}$ and this map induces a group homomorphism $\phi:U(\mcC) \to \widehat{G(\mcC)}$, which is bijective if $\CC$ is modular \cite[Theorem 6.2]{GN}.

\subsubsection{Zesting}

Let $\mcC$ be a modular category and $B\subset G(\mcC)$ a subgroup. Thus, the composition of the restriction map $\widehat{G(\mcC)}\twoheadrightarrow  \widehat{B}$ with the isomorphism $\phi: U(\mcC)\to \widehat{ G(\mcC)}$ defines a $\widehat{B}$-grading of $\mcC$, where $\mcC_0$ is the fusion subcategory  generated by $\{X_i\in \operatorname{Irr}(\mcC): c_{X_i,b}c_{b,X_i}=1, \forall b\in B\}$, that is, $\mcC_0=C_{\mcC}(B)$ the centralizer of $B$ in $\mcC$. Note that
$G(\mcC_0)=\{a\in G(\mcC):c_{a,b}c_{b,a}=1, \forall b\in B\}$. In particular if $B=G(\mcC)$, $A$ is symmetric.

Each $a\in G(\mcC_0)$, defines a $\mcC_0$-bimodule equivalence $L_a:\mcC_\sigma\to \mcC_\sigma, X\mapsto a\otimes X$, with natural isomorphism $c_{a,V}\otimes \id_X: L_a(V\otimes X)\to V\otimes L_a(X)$, for all $\sigma \in \widehat{B}$.

Let $A\subset G(\mcC_0)$ be a subgroup such that the pointed fusion subcategory of $\mcC_0$ generated by  $A$ is symmetric. Thus,  we can assume that the braiding on $A$  is defined by a symmetric bicharacter $c:A\times A\to \{1,-1\}$.

Given $\alpha\in Z^2(\widehat{B},A)$ we  define a new tensor product $\overline{\otimes}_\alpha:\mcC\times \mcC\to \mcC$ as $$\overline{\otimes}_\alpha|_{\mcC_{\sigma}\boxtimes \mcC_{\tau}}=L_{\alpha(\sigma,\tau)}\circ \otimes.$$

By \cite[Proposition 9]{CGPW} the obstruction to the commutativity of the pentagonal diagram of this new tensor product is given by the cohomology class of the following 4-cocycle $O_4(\alpha,c)\in Z^4(\widehat{B},U(1))$,

$$O_4(\alpha,c)(\sigma,\tau,\rho,\eta)=c(\alpha(\sigma,\tau),\alpha(\rho,\eta)),$$ that is,  $O_4(\alpha,c)=\alpha\cup_c\alpha$ (cup product).

Since $\alpha$ is a 2-cocycle, we can assume the innocuous condition \[\alpha(\sigma\tau)\otimes \alpha(\sigma,\tau)=\alpha(\sigma,\tau\rho)\otimes \alpha(\tau,\rho),\] for all $\sigma,\tau,\rho \in \widehat{B}$. Assume that there is $w\in C^3(\widehat{B},U(1))$ such that $\delta(w)=O_4(\alpha,c)$, thus the isomorphisms
$$w_{\sigma,\tau\rho}\id:\alpha(\sigma\tau,\rho)\otimes \alpha(\sigma,\tau)\to \alpha(\sigma,\tau\rho)\otimes \alpha(\tau,\rho)$$ are such that the natural isomorphisms

\begin{equation}\label{new associator}
\hat{a}_{X_\sigma,X_\tau,X_\rho}^w=(\id_{\chi(\sigma,\tau\rho)}\otimes c_{\alpha(\tau,\rho),X_\sigma}\otimes \id_{X_\rho})\circ(w_{\sigma,\tau,\rho}\otimes \id_{X_\sigma\otimes X_\tau\otimes X_\rho}),    
\end{equation}
define an associator with respect to $\bar{\otimes}_\alpha$ and we get  a new $\widehat{B}$-graded fusion category 
$$
\mcC_{(\alpha,w)}:=(\mcC,\bar{\otimes}_\alpha,\hat{a}^w),
$$that we will call a \textit{zesting} of $\mcC$. In case that $\alpha\equiv 1$,  then $w\in Z^3(\widehat{B},U(1))$ is just a 3-cocycle and $\mcC_{(1,w)}$ is called a twisting.

\subsubsection{Zested extensions of a super-modular category}
Let $\mcB$ be a super-modular category and $(\mcC,c)$ a modular closure of $\mcB$. Continuing with the notation of the previous subsection, take $A=B=\{\one, f\}\cong \mathbb{Z}_2$, where $f\in \mcB$ is the
fermion object.  We will identify $A$, $B$ and $\widehat{B}$ with $\mathbb{Z}_2=\{0,1\}$. Let 
$c:A\times A\to \{1,-1\}$ be the non-trivial symmetric
bicharacter, that is $c(f,f)=-1$. Since $H^2(\widehat{B},A)=H^2(\mathbb{Z}_2,\mathbb{Z}_2)\cong \mathbb{Z}_2$, 
\begin{equation}\label{def alpha}
    \alpha(1,1)=f
\end{equation} represents the unique non-trivial  cohomology class. The fourth obstruction in this case is given by the 4-cocycle 
\[O_4(\alpha)(1,1,1,1)=c(\alpha(1,1),\alpha(1,1))=c(f,f)=-1.\]
If we define \[w\in C^3(\widehat{B},U(1)), \  \  \       w(1,1,1)=i,\]  $\delta(w)=O_4(\alpha)$, thus the zesting $\mcC_{(\alpha,w)}$ has associator \begin{equation}\label{associator supermodular zesting}
\hat{a}_{X_\sigma,X_\tau,X_\rho}^b:=w(\sigma,\tau,\rho)\id_{\chi(\sigma,\tau\rho)}\otimes c_{\alpha(\tau,\rho),X_\sigma}\otimes \id_{X_\rho},   
\end{equation}where $\sigma,\tau, \rho\in \widehat{B}$.

\begin{theorem}\label{zested modular closure}
Let $(\CC,c)$ be a modular closure of a super-modular category $(\mcB,f)$. 

\begin{enumerate}

\item Let $\chi: \mathbb{Z}_2\times \mathbb{Z}_2\to \mathbb{Z}_2$ be the non-trivial bicaharcter. With $\alpha$ and $w$ fixed as above, the zesting $\CC_{(\alpha,w)}$ with the natural isomorphism \[c_{X_\sigma,X_\tau}^{\alpha}= (e^{\pi i/4})^{\chi(\sigma,\tau)}\id_{\alpha(\sigma,\tau)}\otimes c_{X_\sigma,X_\tau},\] defines a   modular closure $\CC_\alpha:=(\CC_{(\alpha,w)}, c^\alpha) \supset \mcB$, inequivalent to $(\CC,c) \supset \mcB$. 
\item The $S$ and $T$ matrices of $(\CC_{(\alpha,w)}, c^\alpha)$ are
\begin{equation}\label{S and T matrices of z} \tilde{s}^{\alpha}_{X_\sigma,X_\tau}=i^{\chi(\sigma,\tau)}\tilde{s}_{X_\sigma,X_\tau}, \  \  \  \  \  \theta^{\alpha}_{X_\sigma}=(e^{\pi i/4})^{\chi(\sigma,\sigma)}\theta_{X_\sigma}, \end{equation} for all $X_\sigma\in \mcC_\sigma, X_\tau\in \mcC_\tau,  \sigma,\tau \in \mathbb{Z}_2$.

\item The rule $\CC\mapsto \CC_\alpha$ defines a free action of  $\mathbb{Z}_8$  on the set of equivalence classes of modular closures of $\mcB$. 

\end{enumerate}

\end{theorem}
\begin{proof}
It is straightforward to check the commutativity of the hexagon diagrams. By definition $\CC_\alpha$  is a braided $\mathbb{Z}_2$-extension of $\mcB$. We only need to see the formulas of the new $S$ and $T$ matrices, since they imply that $\CC_\alpha$  is modular. Let $X,Y\in \mcC_{1}$ be defect objects, then $$c^{\alpha}_{Y,X}c^{\alpha}_{X,Y}=(e^{\pi i/4})^2c_{Y,X}c_{X,Y}=ic_{Y,X}c_{X,Y},$$taking the quantum trace we get $\tilde{s}^{\alpha}_{X,Y}=i\tilde{s}_{X,Y}$. Using that for any pre-modular category with $X$ a simple object $\theta_Xd_X=\operatorname{Tr}(c_{X,X})$, we have that $\theta_{X}^\alpha=\operatorname{Tr}(c_{X,X}^{\alpha})/d_X=e^{\pi i/4}\theta_X$ for $X\in\mcC_1$.

It clear from the definition of $\CC_\alpha$ that applying the zesting procedure to $\CC$ eight times returns $\CC$. We only need to check that the action is transitive, which is accomplished by showing the multiplicative central charge $\xi(\CC_\alpha)$ of $\CC_\alpha$ is $e^{\pi i /4}\xi(\CC)$.

By Corollary \ref{gauss sum is zero} the  multiplicative central charge of a modular closure of $\mcB$ is \[\xi(\CC)=\frac{\sum_{X\in \operatorname{Irr}(\CC_1)} \theta_Xd(X)^2}{\sqrt{\operatorname{dim}(\CC)}}.\]
Then \begin{align*}
    \xi(\CC_\alpha)&=\frac{\sum_{X\in \operatorname{Irr}(\CC_1)} \theta_X^\alpha d(X)^2}{\sqrt{\operatorname{dim}(\CC)}}\\
    &= e^{\pi i/4}\frac{\sum_{X\in \operatorname{Irr}(\CC_1)} \theta_X d(X)^2}{\sqrt{\operatorname{dim}(\CC)}}\\
    &=e^{\pi i/4}\xi(\CC).
\end{align*}Since multiplicative central charge is an invariant of pre-modular categories, the elements in the $\mathbb{Z}_8$-orbit of $\CC$  are not equivalent modular closures of $\mcB$.
\end{proof}

\subsection{16-fold way for \texorpdfstring{$PSU(2)_{4m+2}$}{PSU24m+2}}

\begin{theorem} The 16 inequivalent Witt classes of modular closures of the super-modular category $PSU(2)_{4m+2}$ have representatives which can be constructed explicitly.  For $m=0$, there are exactly $16$ modular closures up to ribbon equivalence. 
\end{theorem}

\subsubsection{Modular closures via Theorem \ref{16ways}}
Let  $\CC=SU(2)_{4m+2}$ be the (natural) minimal modular closure of $PSU(2)_{4m+2}$. We first apply the construction of Theorem \ref{16ways}  to $\CC$ to generate $16$ inequivalent minimal modular closures of $PSU(2)_{4m+2}$.  Since the multiplicative central charge of $SU(2)_{4m+2}$ is $e^{3(2m+1)\pi i/(8m+8)}$, the central charges of these minimimal modular closures are $e^{\frac{(6+\nu)m+(3+\nu)}{8m+8}\pi i}$, where $1\leq \nu\leq 16$.

First consider one of the eight Ising theories $\mathcal{I}_j$. We denote the objects by $\1,\sigma,e=\psi$. These $8$ theories are distinguished by $\theta_\sigma=e^{\pi i \nu/8}$ where $\nu=2j+1$ with $0\leq j\leq 7$.  

The associated modular closure $[(\CC\boxtimes\II_j)_{\Z_2}]_0$ of $\mcB=PSU(2)_{4m+2}$ is the trivial component of the $\Z_2$-de-equivariantization of $\CC\boxtimes\II_j$, where the Tannakian category $\EE:=Rep(\Z_2)$ appears as the subcategory generated by $(f,e)$.  By \cite{DGNO1} this is $(\EE^\prime)_{\Z_2}$.  To compute the simple objects of $\EE^\prime$, we look for pairs $(X_i,z)\in\CC\boxtimes\II_j$ so that:

$$\tilde{s}_{(X_i,z),(f,e)}=\tilde{s}_{X_i,f}\tilde{s}_{z,e}=d_id_z.$$

Looking at the respective $S$-matrices we find $\EE^\prime$ has objects:

\begin{enumerate}
 \item $(X_{2i},\1)$, $(X_{2i},e)$ for $0\leq i\leq 2m+1$ and
 \item $(X_{2i+1},\sigma)$ for $0\leq i\leq 2m$.
\end{enumerate}

Now to compute the simple objects in $(\EE^\prime)_{\Z_2}$ we look at the tensor action of $(f,e)$ on $\EE^\prime$.  Under the forgetful functor $F:(\EE^\prime)_{\Z_2}\rightarrow \EE^\prime$ we have:
\begin{enumerate}
 \item $(X_{2i},\1)+(X_{4m+2-2i},e)$ for $0\leq i\leq m$
 \item $(X_{2i},e)+(X_{4m+2-2i},\1)$ for $0\leq i\leq m$
 \item $(X_{2i+1},\sigma)+(X_{4m+2-2i-1},\sigma)$ for $0\leq i\leq (m-1)$ and
 \item $(X_{2m+1},\sigma)$.
\end{enumerate}
The first three types above come from simple objects in $(\EE^\prime)_{\Z_2}$, whereas the last object is the image of a sum of 2 simple objects $Y_1$ and $Y_2$ of equal dimension.  Therefore the rank of $(\EE^\prime)_{\Z_2}$ is $3m+4$. 

The first $2(m+1)$ simple objects in $(\EE^\prime)_{\Z_2}$ coming from $(X_{2i},\1)$ and $(X_{2i},e)$ for simple $X_{2i}\in\mcB$ obviously have dimension $\dim(X_{2i})$, and form the subcategory $[(\EE^\prime)_{\Z_2}]_0\cong\mcB$.   The  $m+2$ simple objects in the odd sector $[(\EE^\prime)_{\Z_2}]_1$ have dimensions $\sqrt{2}\dim(X_{2i+1})$ ($m$ simple objects) and $\frac{\sqrt{2}}{2}\dim(X_{2m+1})$ (2 objects).  

Now let us consider $[(\CC\boxtimes A)_{\Z_2}]_0$ where $A$ is one of the $8$ abelian (pointed) minimal modular closures of $\sVec$.  Explicit realizations of such $A$ can be obtained from (see \cite{RSW}): 1) Deligne products of the rank 2 semion modular category or its complex conjugate (4 theories) 2) the $\Z_4$ modular category and its conjugate 3) the toric code $SO(16)_1$ or 4) the $1$ fermion $\Z_2\times \Z_2$ theory $SO(8)_1$.  We continue to label our chosen fermion by $e$ and the other two non-trivial objects by $a$ and $b$.  In this case a similar calculation gives simple objects in $\EE^\prime$:

\begin{enumerate}
 \item $(X_{2i},\1)$, $(X_{2i},e)$ for $0\leq i\leq 2m+1$ and
 \item $(X_{2i+1},a)$, $(X_{2i+1},b)$ for $0\leq i\leq 2m$.
\end{enumerate}

In this case the tensor action is fixed-point free so we obtain:
\begin{enumerate}
 \item $(X_{2i},\1)+(X_{4m+2-2i},e)$ for $0\leq i\leq m$
 \item $(X_{2i},e)+(X_{4m+2-2i},\1)$ for $0\leq i\leq m$
 \item $(X_{2i+1},a)+(X_{4m+2-2i-1},b)$ for $0\leq i\leq (m-1)$ and $m+1\leq i\leq 2m+1$ and
 \item $(X_{2m+1},a)+(X_{2m+1},b)$.
\end{enumerate}
We see that the rank of $[(\CC\boxtimes A)_{\Z_2}]_0$ is $4m+3$, as expected.

\subsubsection{Explicit data and realizations for modular closures of \texorpdfstring{$PSU(2)_{4m+2}$}{PSU24m+2}}
The 16 minimal modular closures of $PSU(2)_{4m+2}$ can all be constructed from quantum groups.  We record the $S$- and $T$-matrices as they have a fairly simple form.  We group the modular closures into two classes by their ranks: $3m+4$ and $4m+3$.  Notice that for $m=1$ these two cases coincide, so that the constructions below only give $8$ theories: indeed  $SU(2)_6\cong SO(3)_3$.  However, we still obtain $16$ distinct quantum group constructions because $PSU(2)_6$ is equivalent (by a non-trivial outer automorphism) to its complex conjugate: by taking the complex conjugates of each of the $8$ theories constructed (twice) below we obtain a full complement of $16$ modular closures.

The data for the  $8$ modular closures obtained from Ising categories are given in terms of those of the modular category $SO(2m+1)_2$ of rank $3m+4$.   The subcategory $PSO(2m+1)_3$ generated by the objects labeled by integer weights $\lambda\in\Z^{m}$ can be shown to be equivalent to $\overline{PSU(2)}_{4m+2}$ (i.e. the complex conjugate of $PSU(2)_{4m+2}$), with rank $2m+2$.  The other component (with respect to the $\Z_2$ grading) has rank $m+2$ and with simple objects labeled by weights $\mu\in(\frac{1}{2},\ldots,\frac{1}{2})+\Z^m$.

Let $\tilde{S}$ and $\tilde{T}$ be the $S$- and $T$-matrices of $SO(2m+1)_3$, and let $\xi=e^{2\pi i\alpha/8}$ be any $8$th root of unity.  The $8$ rank $3m+4$ minimal modularizations of $PSU(2)_{4m+3}$ have the following data:

$$\tilde{s}_{\lambda,\mu}:=\begin{cases} \tilde{S}_{\lambda,\mu} & \lambda\,\mathrm{or}\,\mu\in\Z^m \\[1em]
\dfrac{\tilde{S}_{\lambda,\mu}}{{\xi}^2} &\lambda,\mu\not\in\Z^m\end{cases}$$

and $$t_{\lambda,\lambda}:=\begin{cases} \tilde{T}_{\lambda,\lambda} & \lambda\in\Z^m\\ \xi \tilde{T}_{\lambda,\lambda} & \lambda\not\in\Z^m.\end{cases}$$

The multiplicative central charges for these theories are $\xi e^{3m(2m+1)\pi i/(8m+8)}$.  
 Although the categories $SO(2m+1)_3$ have been studied (see \cite{gannoncjm}) explicit modular data do not seem to be available.  Direct computation of the data (for example by antisymmetrizations of quantum characters over they corresponding Weyl group) is possible but cumbersome.  For the reader's convenience (and posterity) we provide explicit formulae for $\tilde{S}$ and $\tilde{T}$.  

For a fixed $m$, define $\chi(i,j)=\dfrac{\sin\left(\frac{(i+1)(j+1)\pi}{4m+4}\right)}{\sin\left(\frac{\pi}{4m+4}\right)}$.  Next define the following matrices:

\begin{enumerate}
\item $A_{i,j}:=\chi(2i,2j)$ for $0\leq i,j\leq 2m+1$, so $A$ is $(2m+2)\times (2m+2)$,
\item $B_{k,1}=B_{k,2}=\frac{1}{\sqrt{2}}\chi(2k,2m+1)$ for $0\leq k\leq 2m+1$, so $B$ is $(2m+2)\times 2$,
\item $C_{i,j}:=\sqrt{2}\chi(2i,2j+1)$ for $0\leq i\leq 2m+1$ and $0\leq j\leq m-1$, so $C$ is $(2m+2)\times m$,
\item $D_{i,j}:=\sqrt{\dfrac{m+1}{2}}\dfrac{(-1)^{i+j}}{\sin\left(\frac{\pi}{4m+4}\right)}$ for $0\leq i,j\leq 1$, so $D$ is $2\times 2$.
\end{enumerate}

Now set $$\tilde{S}=\begin{pmatrix}
A & B & C\\
B^T & D &\mathbf{0}\\
C^T & \mathbf{0} & \mathbf{0}
\end{pmatrix}.$$

Define $q=e^{\frac{\pi i}{8m+8}}$.  
The diagonal matrix $\tilde{T}$ has entries:

$$(q^{4(j^2+j)}, 0\leq j\leq 2m+1,q^{-2m^2-m},q^{-2m^2-m},q^{4i^2-(6m^2+9m+4)},1\leq i\leq m).$$

Here the ordering of the simple objects is such that the first $2m+2$ are the objects in $PSO(2m+1)_3\cong PSU(2)_{4m+2}$, i.e. the objects labeled by integral $\mathfrak{so}_{2m+1}$ weights, with corresponding $S$-matrix equal to $A$.   In particular the $2m+2$nd object is the fermion $f$.  The objects corresponding to the columns of $B$ are the two objects in the non-trivial sector that are not fixed under tensoring with the fermion $f$, and the remaining $m$ are each $f$-fixed.

For calibration we point out that for $m=0$ we obtain the Toric Code modular category.

These $8$ categories can be constructed explicitly as follows:

\begin{enumerate}
    \item The construction of $SO(2m+1)_3$ from $U_q\mathfrak{so}_{2m+1}$ with $q=e^{\pi i/(4m+4)}$ depends on a choice of a square root of $q$, and the associativity constraints of each of these can be modified by a $\Z_2$-twist (see \cite{TubaWenzl}) giving the four categories with $\xi^4=1$ above.
    \item By zesting the $4$ theories above (see Section \ref{zest}), we obtain $4$ new non-self-dual categories corresponding $\xi^4=-1$, see Section \ref{zested modular closure}.
\end{enumerate}

Again, let $\xi=e^{2\pi i/8}$ be any $8$th root of unity.  The $8$ rank $4m+3$ minimal modularizations of $PSU(2)_{4m+2}$ have the following data:

$$\tilde{s}_{i,j}:=\begin{cases} \dfrac{\sin\left(\frac{(i+1)(j+1)\pi}{4m+4}\right)}{\sin\left(\frac{\pi}{4m+4}\right)} & 2\mid ij,\\[1em]
\dfrac{\sin\left(\frac{(i+1)(j+1)\pi}{4m+4}\right)}{{\xi}^2\sin\left(\frac{\pi}{4m+4}\right)} & 2\nmid ij \end{cases}$$

and
$$t_{j,j}:=\begin{cases}
e^{\frac{\pi i(j^2+2j)}{8m+8}} & 2\mid j,\\
\xi e^{\frac{\pi i(j^2+2j)}{8m+8}} & 2\nmid j.
\end{cases}$$

The multiplicative central charges for these theories are $\xi e^{3(2m+1)\pi i/(8m+8)}$.
These categories can be realized as follows:
\begin{enumerate}
    \item  $SU(2)_{4m+2}$ is obtained from $U_q\mathfrak{sl}_2$ with $q=e^{\pi i/(4m+4)}$ by choosing the square root of $q$ with the smallest positive angle with the $x$-axis.  The other choice provides a distinct category.  The associativity constraints of these categories can be twisted in two ways using \cite{KW} to obtain a total of 4 categories.  These correspond to $\xi^4=1$.
    \item By zesting the  4 theories above (see Section \ref{zest}) we obtain the 4 non-self-dual modular categories, corresponding to $\xi^4=-1$, cf. Section \ref{zested modular closure}. Alterernatively, we can use the results of \cite[Theorem 5.1]{OS} to see that $PSU(2)_{4m+2}$ and the ``mirror" category to $PSU(4m+2)_2$ are equivalent as ribbon categories.  Since $SU(4m+2)_2$ is obviously a minimal modular extension of $PSU(4m+2)_2$ we can proceed as above to find 4 distinct versions: two for the choice of a (square) root of $q$ and another two from the two Kazhdan-Wenzl twists that preserve $PSU(4m+2)_2$.
\end{enumerate}

\section{A graphical calculus for zesting}\label{zest}
\newcommand{\figref}[3]{
  \begin{figure}
    \begin{center}
      \includegraphics[#3]{#1}
      \caption{\label{#1}#2}
    \end{center}
  \end{figure}
}

In this section, given a supermodular category $\mathcal B$ with modular closure $\mathcal C$, we construct seven other modular closures using the graphical calculus for $\mathcal C$. Another, more general, approach would be to apply results of \cite{LKW} and Definition/Proposition 2.15 in \cite{M3} directly to compute categorical data for all sixteen modular closures. That approach, however, requires explicit computation of idempotent completions; the approach considered here provides computational simplicity at the cost of some generality.

Let $\mathcal C$ be a $\mathbb Z_2$-graded unitary modular category over $\mathbb C$, with Grothendieck semiring $R$, containing a pointed object $e$ of order two in $\mcC'$.

The object $e$ generates a subcategory equivalent as a braided fusion category to $\operatorname{Rep}(\mathbb Z_2$) or $\sVec$. Since $\operatorname{dim}(e) = 1$, we have
$c_{e,e} = \theta_e Id_{e \otimes e},$
with $\theta_e = \pm 1$ depending on if $e$ is a boson or fermion.

Let $\mathcal C_0$ and $\mathcal C_1$ denote the trivial and nontrivial gradings of $\mathcal C$ respectively. An object or morphism is {\em even} (resp. {\em odd}) if it lies in $\mathcal C_0$ (resp. $\mathcal C_1$). Every object $x \in Ob(\mathcal C)$ is (isomorphic to) a direct sum of even and odd objects. Given two such even-odd direct sum decompositions $x = x_0 \oplus x_1$ and $y = y_0 \oplus y_1$, every $f:x\to y$ decomposes uniquely as $f = f_0 \oplus f_1$, where $f_0:x_0 \to y_0$ and $f_1:x_1 \to y_1$.

\subsection{Zested fusion rules}

There is a bifunctor of categories $\boxtimes:\mathcal C \times \mathcal C \to \mathcal C$ which acts on simple objects $x_1, x_2 \in \operatorname{Ob}(\mathcal C)$ as follows:
\[x_1 \boxtimes x_2 :=
\left\{
\begin{array}{ll}
  x_1 \otimes x_2 & \mbox{ if at least one of }x_1, x_2\mbox{ lies in }\operatorname{Ob}(\mathcal C_0), \\  
  (x_1 \otimes e) \otimes x_2 & \mbox{ otherwise.}
\end{array}
\right.
\]

The operation of $\boxtimes$ on even and odd morphisms is defined by

\[f_1 \boxtimes f_2 :=
\left\{
\begin{array}{ll}
  f_1 \otimes f_2 & \mbox{ if at least one of }f_1, f_2 \mbox{ is even}, \\
  (f_1 \otimes Id_e) \otimes f_2 & \mbox{ if both }f_1\mbox{ and } f_2 \mbox{ are odd}.
\end{array}
\right.
\]

The functor $\boxtimes$ gives (isomorphism classes of) objects in $\mathcal C$ a $\mathbb Z^+$-based semiring structure $R^\boxtimes$.

It is convenient to distinguish instances of $e$ which are introduced by the $\boxtimes$ operator from other instances by referring to them as {\em gluing objects}.

\subsection{Associativity}
Let $\alpha$ be the associator of $\mathcal C$, $\lambda$ and $\rho$ the triangle isomorphisms, and $c$ the braiding.

Fix two constants $l,r \in \mathbb C$. For each triple of simple objects $a,b,c \in Ob(\mathcal C)$, define the map $\beta_{a,b,c}:(a \boxtimes b) \boxtimes c \to a \boxtimes (b \boxtimes c)$ as follows. Note that here and in the rest of this section we use the composition of arrows convention, so that $f \circ g$ has domain $\operatorname{dom}(f \circ g) = \operatorname{dom}(f)$.
\begin{itemize}
\item If at most one of $a,b,c$ is odd, $\beta_{a,b,c} = \alpha_{a,b,c}$.
\item If $c$ alone is even, $\beta_{a,b,c} = \alpha_{a\otimes e,b,c}$.
\item If $a$ alone is even, \[
\beta_{a,b,c} = ((a \otimes b) \otimes e) \otimes c \xrightarrow{\alpha_{a \otimes b, e, c}} (a \otimes b) \otimes (e \otimes c) \xrightarrow{ \alpha_{a, b, e \otimes c}}\] \[ a \otimes (b \otimes (e \otimes c))\xrightarrow{Id_a \otimes \alpha^{-1}_{b,e,c}} a \otimes ((b \otimes e) \otimes c).\]
\item If $b$ alone is even, 
\[\beta_{a,b,c} = ((a \otimes b) \otimes e) \otimes c \xrightarrow{l (\alpha_{a,b,e} \otimes Id_c)} (a \otimes (b \otimes e)) \otimes c \xrightarrow{(Id_a \otimes c^{-1}_{e,b}) \otimes Id_c}\] \[ (a \otimes (e \otimes b)) \otimes c \xrightarrow{\alpha^{-1}_{a,e,b} \otimes Id_c}((a \otimes e) \otimes b) \otimes c \xrightarrow{ \alpha_{a\otimes e, b,c}} (a \otimes e) \otimes (b \otimes c).\]
\item if $a,b,c$ are all odd, 
\[\beta_{a,b,c} = ((a \otimes e) \otimes b) \otimes c \xrightarrow{r (\alpha_{a,e,b} \otimes Id_c)} (a \otimes (e \otimes b)) \otimes c \xrightarrow{ \alpha_{a,e \otimes b, c}}\] \[ a \otimes ((e \otimes b) \otimes c) \xrightarrow{Id_a \otimes (c_{e,b} \otimes Id_c)} a \otimes ((b \otimes e) \otimes c).\]
\end{itemize}
One may interpret the definition pictorially by applying a factor of $r$ (resp. $r^{-1}$) whenever a gluing object is slid to the right (resp. left) over an odd object due to reassociation.

Extend these definitions to all triples of objects via direct sum decompositions.

Then $(\mathcal C, \boxtimes, \beta, \lambda, \rho)$ is a monoidal category, if $l$ and $r$ are nonzero, $\beta$ is natural with respect to morphisms, $\lambda$ and $\rho$ are natural isomorphisms $\lambda_x: \triv \otimes x \to x$ and $\rho_x:x \otimes \triv \to x$ satisfying the triangle axioms (it is well-known that such morphisms always exist if the other conditions in the definition can be satisfied), and for all $a,b,c,d \in \mathcal C$, the following coherence property holds:
\begin{equation}\label{coherence-beta}\beta_{a \otimes b,c,d} \circ \beta_{a, b, c\otimes d} = (\beta_{a,b,c} \otimes Id_d) \circ \beta_{a,b \otimes c, d} \circ (Id_a \otimes \beta_{b,c,d}).
\end{equation}

Naturality of $\beta$ with respect to morphisms $f:a \to b$ follows from naturality of associativity $\alpha$ and $c$ with respect to morphisms; the constants on either side of the naturality equation cancel by a parity argument. Furthermore, by the coherence property and naturality of the braiding $c$ over $\alpha$, the validity of each instance of Equation~\ref{coherence-beta} is determined entirely by the following:
\begin{itemize}
\item The values of $l$ and $r$,
\item The domain and range (equal on both sides of each equation),
\item In the case of four odd objects, the braiding of the two gluing objects.
\end{itemize}

The powers of $l$ and $r$ which occur on each side of Equation~\ref{coherence-beta}, as well as the number of instances of $c_{e,e}$, depend only on the parity of the objects. If not all of $a,b,c,d$ are odd, the only possible relation on $r$ and $l$ is that $l^2 = l$, obtained in the odd-even-even-odd case. Thus we set $l = 1$.

If all of $a,b,c,d$ are odd, then
\[((a \boxtimes b) \boxtimes c) \boxtimes d = ((((a\otimes e) \otimes b)\otimes c) \otimes e) \otimes d.\]
In this case, the right hand side of the coherence equation differs from the left in that it has a factor of $r^2$ and an exchange $c_{e,e}$ of the two gluing objects. Since $c_{e,e} = \theta_e Id_{e \otimes e}$, we obtain the following:

\begin{lemma}, Let $\mcC^\boxtimes = (Ob(\mathcal C), \boxtimes, \beta, \lambda, \rho)$. When $l=1$ and $r^2 = \theta_e$, $\mcC^\boxtimes$ is a monoidal category.
\end{lemma}

Note that there is a canonical isomorphism $R \cong R^{\boxtimes \boxtimes}$. If $\mcC$ is skeletal, then $\mcC$ and $\mcC^{\boxtimes \boxtimes}$ have isomorphic Grothendieck rings and identical associators, except that the odd-odd-odd associators in $\mathcal{C}^{\boxtimes \boxtimes}$ differ from those in $\mathcal C$ by a factor of $r^2 \theta^2_e = \theta_e$. One consequence is that if $\theta_e = 1$, then $\mcC \cong \mcC^{\boxtimes \boxtimes}$, and otherwise applying the construction twice gives $\mathcal C \cong \mathcal C^{\boxtimes \boxtimes \boxtimes \boxtimes}$. Furthermore, when $\theta_e = -1$, one finds that $\mathcal{C}^{\boxtimes \boxtimes \boxtimes}$ is equivalent to what would result from $\mathcal{C}^{\boxtimes}$ if the other choice of sign for $r$ were made.

It is less clear whether or not additional equivalences exist. Ultimately, we will obtain eight modular categories, the non-equivalence of which is shown by the central charge, and at the level of fusion categories, no such invariant exists. As our interest is in the modular structure we do not attempt to completely specify equivalences at the level of fusion categories.

\subsection{Rigidity}

Let $x \in \operatorname{Ob}(\mathcal C)$ be simple. Let
\[x^\boxast =\left\{
\begin{array}{ll}x^* & \mbox{ if }x\mbox{ is even,} \\
  e^* \otimes x^* & \mbox{ otherwise.}
\end{array}
\right.
\]

Define the maps $\operatorname{ev}_x^\boxtimes : x \boxtimes x^\boxast \to \triv$ and $\operatorname{coev}_x^\boxtimes : \triv \to x^\boxast \boxtimes x$ such that if $x$ is even we have $\operatorname{coev}_x^\boxtimes = \operatorname{coev}_x$ and $\operatorname{ev}_x^\boxtimes = \operatorname{ev}_x$, and if $x$ is odd,
\[\operatorname{ev}_x^\boxtimes = (\Id_x \otimes \alpha^{-1}_{e,e^*,x^*}) \circ (\Id_x \otimes (\operatorname{ev}_e \otimes \Id_{x^*})) \circ (\Id_x \otimes \lambda(x^*)) \circ \operatorname{ev}_x,\]

\[\operatorname{coev}_x^\boxtimes = r^{-1} \operatorname{coev}_x \circ (\rho^{-1}(x^*) \otimes \Id_x) \circ ((\Id_{x^*} \otimes \operatorname{coev}_e) \otimes \Id_x)\]
\[\circ (\alpha^{-1}_{x^*,e^*,e} \otimes \Id_x) \circ ((c^{-1}_{e^*,x^*} \otimes \Id_e) \otimes \Id_x.\]

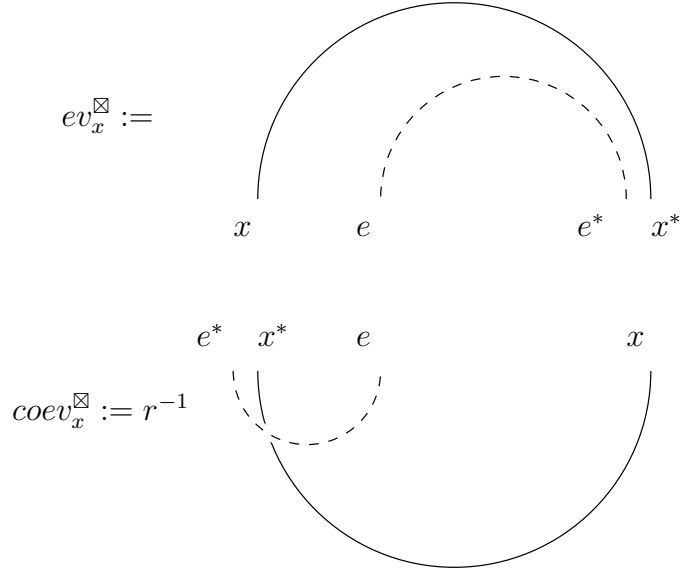
\figref{fig-birth-death}{The birth and death on odd $x$ in $\mathcal C^\boxtimes$.}{height=3in}

See Figure~\ref{fig-birth-death}. Factors of $r$ again algebraically count the crossings of gluing strands over odd strands. This feature will persist throughout the construction.

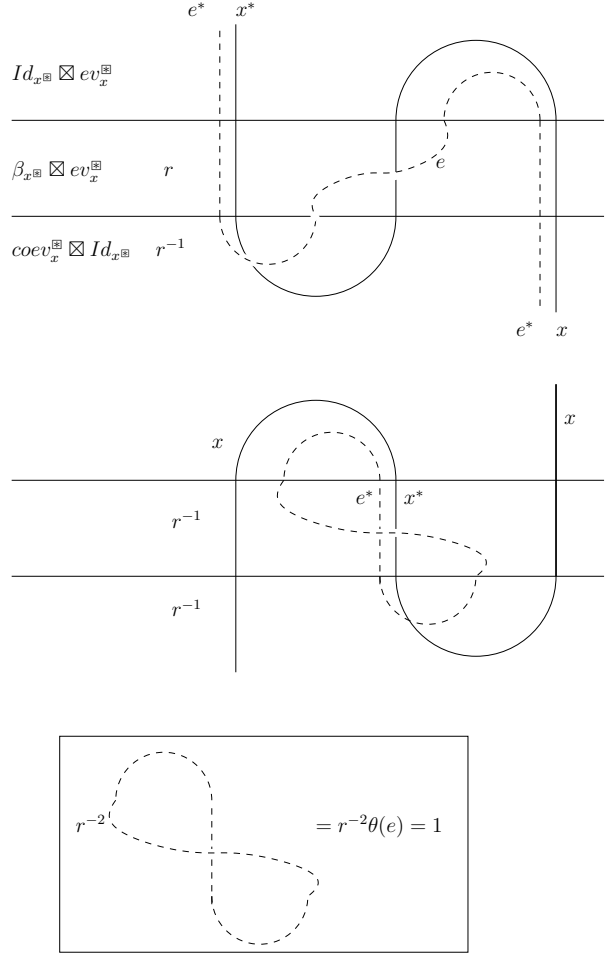
\figref{fig-rigid}{Figures for the rigidity equations.}{height=5in}

We have
\[\rho^{-1}_x \circ (\Id_x \otimes \operatorname{coev}^\boxtimes_x) \circ \beta^{-1}_{x,x^\boxast,x} \circ (\operatorname{ev}_x \otimes \Id_x) \circ \lambda_x = \Id_x\]
by standard graphical calculus techniques, since the morphism $\operatorname{coev}_e \circ c_{e^*,e} \circ \operatorname{ev}_e$ evaluates to $\theta_e$ and there is a factor of $r^{-1}$ from $\beta^{-1}_{x,x^\boxast,x}$. See Figure~\ref{fig-rigid}.

Along similar lines,
\[\lambda^{-1}_{x^\boxast} \circ (\operatorname{coev}_x^\boxtimes \otimes \Id_{x^\boxast}) \circ \beta_{x^\boxast,x,x^\boxast} \circ (\Id_{x^\boxast} \otimes \operatorname{ev}_x^\boxtimes) \circ \rho_{x^\boxast} = \Id_{x^\boxast},\]
since the factor of $r$ in $\beta_{x^\boxast,x,x^\boxast}$ cancels the constant in $\operatorname{coev}_x^\boxtimes$. See Figure~\ref{fig-rigid}.

Thus $\mathcal C^\boxtimes$ is rigid.

Clearly $\mathcal C^\boxtimes$ is a fusion category with fusion subcategory $(\mathcal C_0, \otimes|_{\mathcal C_0}, \alpha|_{\mathcal C_0}, \lambda|_{\mathcal C_0}, \rho|_{\mathcal C_0})$.

\subsection{Graphical Calculus}

Let $f$ be a composition of identity-tensored reassociations $\beta$ on a product $x_1\otimes \cdots\otimes x_n$ of even or odd objects $x_i$. In terms of $\mathcal C$, $f$ is some power $r^k$ of $r$ times a composition of identity-tensored maps $\alpha$ and instances of $c$. In the strict picture calculus for $\mathcal C$, $f$ is represented, up to factor $r^k$, by a braiding of the $n$ tensored objects $x_i \in \mathcal C$  with at most $\lfloor \frac n 2 \rfloor$ gluing objects. The braiding satisfies the following properties:

\begin{enumerate}
\item The $x_i$ braid trivially with each other.
\item \label{interpose-property}At each stage of the composition, (before or after an instance of $\beta$), each pair of gluing objects is separated by an odd object $x_i$.
\item The number of gluing objects is always half the number of odd $x_i$, rounded down.
\end{enumerate}

The following proposition asserts that any picture satisfying the above properties represents a well-defined morphism in $\mathcal C^\boxtimes$.
\begin{proposition} Let $X_e$ be a multiset of even objects in $\mathcal C$, and $X_o$ a multiset of odd objects in $\mathcal C$. Let $x_1 o_1 x_2 \ldots o_{n-1} x_n$ be a formal string, with $n\ge 2$, satisfying the following conditions:
  \begin{enumerate}
  \item $\{x_1, \ldots, x_n\}$ is the multiset union of $X_e$ and $X_o$,
  \item each $o_i$ is either $\otimes$ or $\boxtimes$,
  \item $\boxtimes$ appears $\lfloor \frac{|X_o|}{2} \rfloor$ times,
  \item If $j > i$ and $o_i = o_j = \boxtimes$, then for some $i < k \le j$, $x_k \in X_o$.
  \end{enumerate}

  Then the following hold:
  \begin{enumerate}
  \item \label{assoc_exists}There exists an association of the operators such that $o_i = \boxtimes$ iff both arguments of $o_i$ are odd.
  \item \label{assoc_parity}Any sequence of (identity-tensored) $\beta$ instances connecting two such associations consists of a sequence of maps $\alpha$ and braidings of the gluing objects over the $x_i$, multiplied by $\theta_e^k$, where $k$ is the sign of the permutation of the gluing objects among themselves.
  \item \label{assoc_trivial}Any two such associations are connected by a sequence of $\beta$ instances which trivially permute the gluing objects. 
  \end{enumerate}
\end{proposition}

\begin{proof}
  First, suppose that $X_e$ is empty.

  By a simple counting argument, there is a pair $(o_i, o_{i+1})$ such that exactly one of $o_i$ and $o_{i+1}$ is $\boxtimes$. Associate to obtain $(x_i \boxtimes x_{i+1}) \otimes x_{i+2}$ or $x_i \otimes (x_{i+1} \boxtimes x_{i+2})$, which is odd in either case, and induct on $n$. This proves (\ref{assoc_exists}) when $X_e$ is empty.

  If $X_e$ is not empty, partially associate the string so that it forms a product of maximal substrings $s_j$ subject to the following conditions:
  \begin{enumerate}
  \item No $s_j$ contains $\boxtimes$,
  \item Each $s_j$ contains exactly one element of $X_o$, with multiplicity.
  \end{enumerate}
  Tensor products within each $s_j$ involve $\otimes$ only, and one may reduce to the previous case. This proves (\ref{assoc_exists}).

  The braiding induced in the picture calculus for $\mathcal C$ is trivial unless there is a reassociation $\beta_{a,b,c}$, where $a$ and $c$ are odd and $b$ contains two or more elements, counted with multiplicity, of $X_o$. Then $b = b_1 \otimes b_2$ for some $b_1$ and $b_2$. In the picture calculus for $\mathcal C$, $\beta_{a,b,c}$ moves the gluing object over the strands of $b$, rightward if $b$ is odd and leftward if $b$ is even.

  By associativity, one may replace $\beta_{a,b,c}$ with
  \[(a \boxtimes (b_1 \boxtimes b_2)) \boxtimes c \xrightarrow{\beta^{-1}_{a,b_1,b_2} \boxtimes \Id_c} ((a \boxtimes b_1) \boxtimes b_2) \boxtimes c \xrightarrow{\beta_{a \boxtimes b_1, b_2, c}} (a \boxtimes b_1) \boxtimes (b_2 \boxtimes c) \xrightarrow{\beta_{a, b_1, b_2 \boxtimes c}}\] \[ a \boxtimes (b_1 \boxtimes (b_2 \boxtimes c)) \xrightarrow{\Id_a \otimes \beta^{-1}_{b_1,b_2,c}} a \boxtimes ((b_1 \boxtimes b_2) \boxtimes c).\]

  In terms of the picture calculus for $\mathcal C$, this has the following effects. If $b_1$ and $b_2$ are not both odd, the braiding of the gluing object over the strands of $b$ is replaced by braidings in the same direction over $b_1$ and $b_2$ individually, and the power of $r$ is not changed. If $b_1$ and $b_2$ are both odd, the rightward braiding of the gluing object over $b_1 \boxtimes b_2 = b_1 \otimes e \otimes b_2$ is replaced with rightward braidings over $b_1$ and $b_2$, along with a factor $r^2 = \theta_e$. The new picture calculus diagram differs topologically from the old in that a single crossing of gluing objects has been replaced by $\Id{e \otimes e}$.

  Repeating this process until one obtains a sequence of identity tensored maps $\beta_{a_i,b_i,c_i}$ such that each $b_i$ contains at most one element of $X_o$, one obtains (\ref{assoc_parity}) and (\ref{assoc_trivial}).
\end{proof}

Notes:
\begin{itemize}
\item By the penultimate paragraph of the previous proof, in the $\mathcal C$-picture calculus, each braiding of a gluing object over an odd strand may be assumed to result from a single odd-odd-odd instance of $\beta$. A morphism in $\mcC^\boxtimes$ inherits, for each such braiding, a factor of $r$ or $r^{-1}$ when the braiding is $c_{e,x}$ or $c_{e,x}^{-1}$ respectively. Thus one may represent reassociativity morphisms in $\mathcal C^\boxtimes$ in the (strict) picture calculus for $\mathcal C$ by adopting the convention that for each $c_{e,x}$ involving a gluing object one multiplies by a factor of $r$ and inversely. Under this convention, any two reassociations with the same picture calculus representations for the domain and codomain become equal.

\item If a tensored object $x_i$ happens to be isomorphic to $e$, but is not introduced as part of an instance of $\boxtimes$, it does not induce a factor of $r$ when it braids with odd objects.

\item We have not shown that there is always a sequence of reassociations in which odd-odd-odd instances of $\beta$ do not occur. Underlying reassociations in $\mathcal C^\boxtimes$ may move the gluing objects. However, there is a way to do it such that the resulting braiding is trivial, and in this case the factors of $r$ all cancel.

\item The braiding of gluing objects with elements of $\mathcal C$ is not natural with respect to picture morphisms. If $x$ and $y$ are strict (i.e. formal) tensor products of even and odd objects, and $f:x \to y$ is a picture morphism such that $x \times y$ has $2$ mod $4$ odd strands, then $c_{e,x} \circ f = - f \circ c_{e,y}$ by a crossing counting argument. For this reason, gluing objects must be distinguished from non-gluing instances of the same object.
\end{itemize}

\subsection{Pivotal and Spherical structure}

Let $\phi$ be the pivotal structure on $\mathcal C$. For any object $x \in \operatorname{Ob}(\mathcal C^{\boxtimes})$, we have
\[x^{\boxast\boxast} =
\left\{
\begin{array}{ll}
  x^{**} &\mbox{if }x\mbox{ is even,} \\
  e^* \otimes (x^{**} \otimes e^{**})&\mbox{if }x\mbox{ is odd.}
\end{array}
\right.
\]

The above-defined rigidity structure on  $\mcC^\boxtimes$ defines a dual functor $\boxast$. We show pivotality using the picture calculus as follows:

Let $f:a \boxtimes b \to c$ be a morphism, with $a$ and $b$ odd. Thus $c$ is even. One may compute the "picture double dual" $\tilde f^{\boxast\boxast}$ of  $f$ (or, similarly, any fusion-category-level picture morphism) as follows:

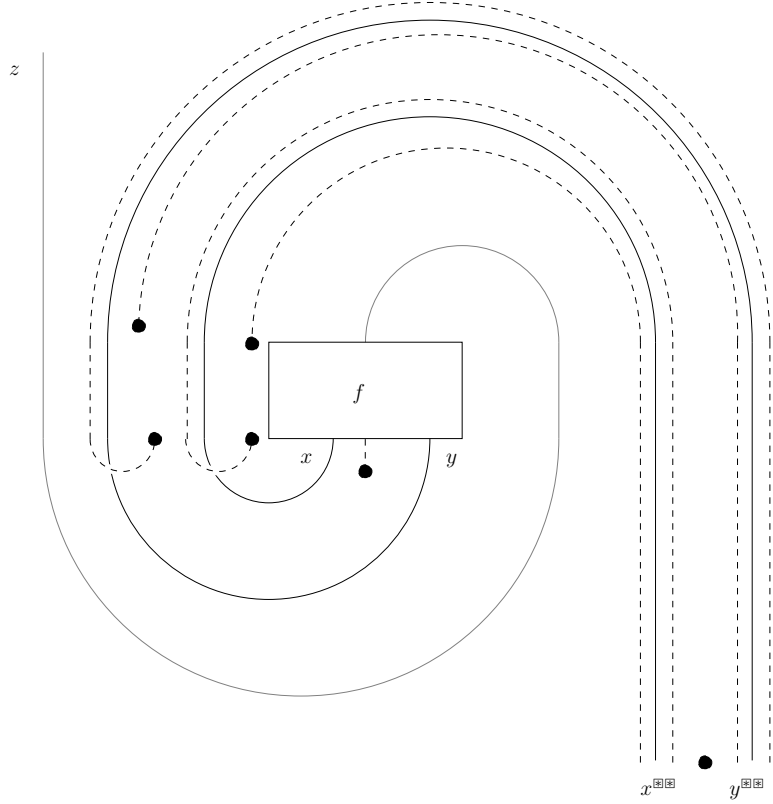
\figref{fig-double-dual-1}{An incomplete picture morphism. The gluing strands terminating at the black dots need to be connected to other strands in order to define a composition. Different compositions may result, but any two give the same morphism up to crossing factors.}{width=4in}
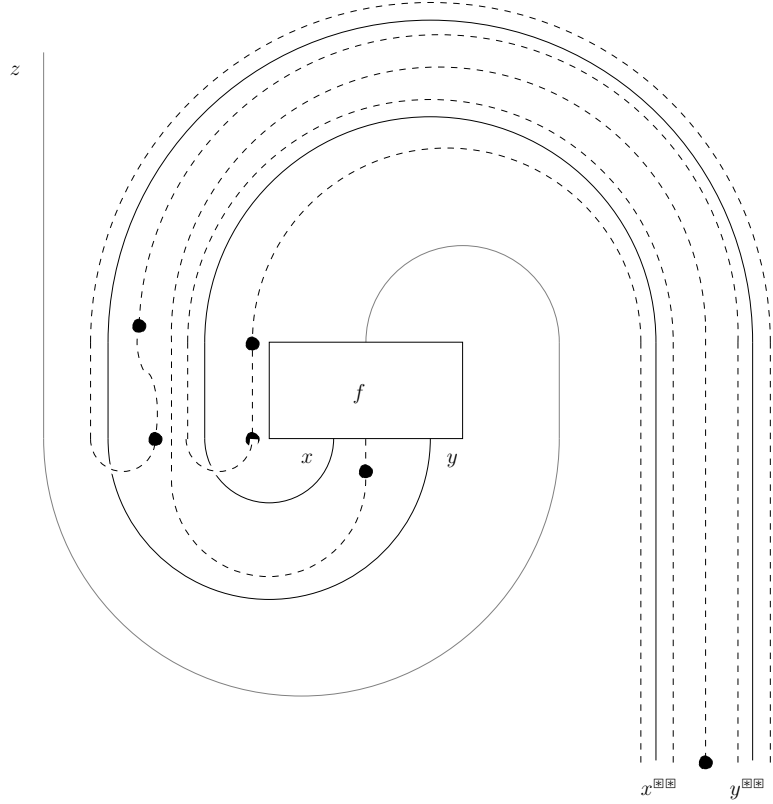
\figref{fig-double-dual-2}{One way to connect the gluing objects. The constant factor is $r^{-2} = \theta_e$ since a gluing object crosses an odd object in each of $\operatorname{ev}^{\boxtimes}_{\operatorname{coev}^\boxtimes}$ and in $\operatorname{ev}^{\boxtimes}_{a^\boxtimes}$. If you don't like the presence of births, deaths, and pivotal isomorphisms on the gluing objects, connect the gluing objects for the domains of $f^{**}$ and $f^{\boxast\boxast}$ along a straight line path, and verify that after accounting for constant factors the same morphism results.}{width=4in}
\begin{enumerate}
\item Draw the usual picture double dual morphism, ignoring gluing objects except as they appear in births, deaths, the domain of $f^{\boxast\boxast}$, and the domain of $f$. See Figure~\ref{fig-double-dual-1}.
\item Connect the gluing objects in any way desired, consistent with the positioning rules. See Figure~\ref{fig-double-dual-2} and its caption for an example.
\item Apply the crossing rules to obtain the appropriate constant factor. In the case of Figure~\ref{fig-double-dual-2}, the factor is $r^{-2} = \theta_e$. 
\end{enumerate}

For each simple object $x$, define $\phi^\boxtimes_x:x^{\boxast \boxast} \to x$ such that
\[\phi^\boxtimes_x =
\left\{
\begin{array}{ll}
  \phi_x & \mbox{if }x\mbox{ is even,}\\
  \begin{array}{l}r^{-1}(\Id_{e^*} \otimes c^{-1}_{e^{**},x^{**}}) \circ \alpha^{-1}_{e^*,e^{**},x^{**}} \\ \circ (\operatorname{ev}_{e^*} \otimes \Id_{x^{**}})\circ \lambda_{x^{**}} \circ \phi_x\end{array} & \mbox{if }x\mbox{ is odd.}
\end{array}
\right.
\]

For odd $x$, the inverse of this map is
\[(\phi^\boxtimes)^{-1}_x = \begin{array}{l}r \phi_x^{-1} \circ \lambda^{-1}_{x^{**}} \circ (\operatorname{coev}_e \otimes \Id_{x^{**}}) \circ \alpha_{e^*,e,x^{**}}\\ \circ (\Id_{e^*} \otimes (\phi^{-1}_e \otimes \Id_{x^{**}})) \circ (\Id_{e^*} \otimes c_{e^{**},x^{**}})\end{array}.\]

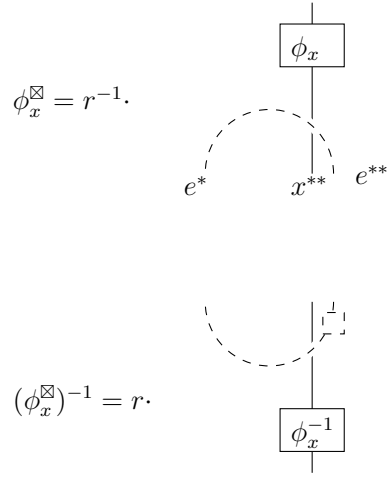
\figref{fig-pivotal}{The pivotal structure for odd $x$ in $\mathcal C^\boxtimes$.}{height=2.5in}

See Figure~\ref{fig-pivotal}.

In a fusion category, the double dual functor $F$ is always isomorphic to the identity as a non-monoidal functor (in a skeletal category, rigidity and semisimplicity imply that the double dual is the identity on the nose). In this case, for any morphism $f:x \to y$, it is clear that
\[f = (\phi^\boxtimes)^{-1}_x \circ f^{\boxtimes \boxtimes} \circ \phi^\boxtimes_y\]
by standard picture calculus techniques (in particular pivotal structure properties of $\phi$ in $\mcC$ and removing loops).

It remains to show that $\phi$ satisfies the monoidal condition:
\[(\phi^\boxtimes)^{-1}_a \otimes (\phi^\boxtimes)^{-1}_b \circ F_2(a \otimes b) \circ \phi^\boxtimes_{a \otimes b} = Id_{a \otimes b}\]

This is done in the usual picture calculus way: Let $c = a \otimes b$, $g:a \otimes b \to c$, $g = Id_{a \otimes b}$, and let $F$ be the double dual functor on $\mcC^\boxtimes$. It is easy to verify that $F_2(a,b) = \tilde g^{\boxast \boxast}$. Breaking up $c$ into its simple object decomposition and applying compatibility of direct sum with tensor product, one has that $\phi^\boxtimes$ is a pivotal category if for all objects $a$ and $b$, simple objects $c$, and morphisms $f:a \otimes b \to c$, we have the following:

\[((\phi^\boxtimes_a)^{-1} \boxtimes (\phi^\boxtimes_b)^{-1}) \circ \tilde f^{\boxast \boxast} \circ \phi^\boxtimes_c = f.\]
This again holds by picture calculus techniques: the case where $a$,$b$,and $c$ are all even follows directly by pivotality in $\mathcal C$, and the case where $a$ and $b$ have opposite parity follows by arguments similar to the above.

Thus the maps $\phi^\boxtimes$ give $\mathcal C^\boxtimes$ a pivotal structure.

Figure \ref{fig-spherical} shows that under this structure, the left and right quantum dimensions of odd objects $x$ in $\mathcal C^\boxtimes$ are equal to the corresponding dimensions in $\mathcal C$. Thus $\mathcal C^\boxtimes$ is a spherical category with $\phi^\boxtimes$ a spherical pivotal structure.

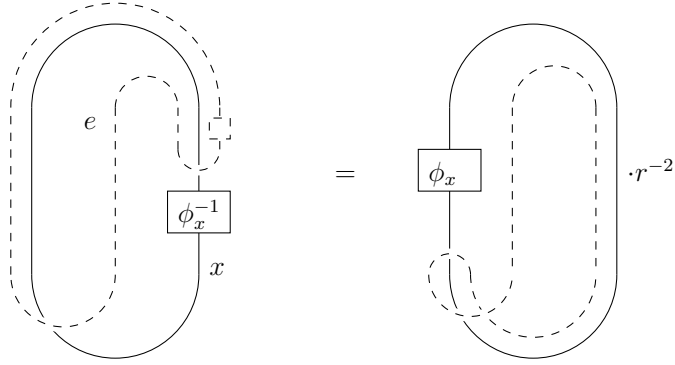
\figref{fig-spherical}{The left and right quantum dimensions in $\mathcal C^\boxtimes$ are equal.}{width=3.5in}
\subsection{Braiding}

For this section we will need some information from the unitary and modular structure of $\mathcal C$. Additionally, we now assume $\mathcal C$ is the modular closure of a supermodular category, and thus $\theta_e = -1$.

\begin{lemma}
  Let $x$ be an odd object in $\mathcal C$. Then $\tilde{s}_{e,x} = -d_x$.
\end{lemma}
\begin{proof}
  For any simple object $y$, 
  \[\frac{(\tilde{s}_{e,y})^2}{d_y} = \tilde{s}_{1,y} = d_y.\]
  Thus we must have
  \[\tilde{s}_{e,y} = \pm d_y.\]

  By assumption, if $y$ is even, the braiding is symmetric, and $\tilde{s}_{e,y} = 1$. In order for $\mathcal C$ to be modular, there must be at least one odd simple object $x_0$ such that $\tilde{s}_{e,x_0} = -d_{x_0}$. But then
  \[\frac{\tilde{s}_{e,x}\tilde{s}_{e,x_0}}{d_e} = \sum_c N_{x,x_0}^c \tilde{s}_{e,c}.\]
  Since $\mathcal C$ is unitary and $N_{x,x_0}^c$ is nonzero only when $c$ is even, in which case $\tilde{s}_{e,c} = d_c$, $\tilde{s}_{e,x}$ must be negative.
\end{proof}

For each pair of objects $x, y$ in $\mathcal C$, and constant $b$, Define $c^\boxtimes_{x,y}:x \boxtimes y \to y \boxtimes x$ such that
\[c^\boxtimes_{x,y} =
\left\{
\begin{array}{ll}
  c_{x,y} & \mbox{ if at least one of }x\mbox{ or } y \mbox{ is even,} \\
  b (c_{x,e} \otimes \Id_y) \circ c_{e\otimes x,y} \circ \alpha^{-1}_{y,e,x} & \mbox{ otherwise.}
\end{array}
\right.
\]

Then $c^\boxtimes$ gives a braiding iff it is natural and satisfies the hexagon equations. Naturality follows by semisimplicity since $c^\boxtimes$ is an isomorphism and is compatible with direct sums, properties it inherits from $c$. The hexagon equations hold if and only if the following two conditions hold for all simple objects $x,y,z,w$ and morphisms $f:x \boxtimes y \to w$,:
\begin{enumerate}
\item $\beta_{x,y,z}^{-1} \circ (f \boxtimes \Id_z) \circ c^\boxtimes_{w,z} = (\Id_x \boxtimes c^\boxtimes_{y,z}) \circ \beta^{-1}_{x,z,y} \circ (c^\boxtimes_{x,z} \boxtimes \Id_y) \circ \beta_{z,x,y} \circ (\Id_z \boxtimes f)$,
\item $\beta_{z,x,y} \circ (\Id_z \boxtimes f) \circ c^\boxtimes_{z,w} = (c^\boxtimes_{z,x} \boxtimes \Id_y) \circ \beta_{x,z,y} \circ (\Id_x \boxtimes c^\boxtimes_{z,y}) \circ \beta_{x,y,z}^{-1} \circ(f \boxtimes \Id_z)$.
\end{enumerate}

Writing out the definitions in terms of $\otimes$, $\alpha$ and $c$, one finds that if at least one of $x$, $y$ or $z$ is even, these equations both follow from naturality properties in the orginal category and cancelling factors $b$.

If $x,y$ and $z$ are all odd, in the first equation, after applying picture calculus operations one obtains $r = b^2$, so we must have $b$ a square root of $r$. In the second equation, we obtain $r^{-1}$ on the left hand side, $b^2$ again on the right hand side, and the morphisms differ by a full twist of the gluing object around $z$. Since $\operatorname{Hom}(e \otimes z, e \otimes z)$ is one dimensional,
\[c_{e,z} \circ c_{z,e} = \frac{\tilde{s}_{e,z}}{d_ed_z} \Id_{e \otimes z} = -\Id_{e \otimes z}.\]
Thus the second equation holds iff
\[r^{-1} = -b^2.\]
Since $b^2 = r$ and $r^2 = \theta_e = -1$, the braid equations are satisfied.

\subsection{\texorpdfstring{$S$-}{S} and \texorpdfstring{$T$-}{T} Matrices}\label{stzested}
Here we describe the $S$- and $T$-matrices for $\mathcal C^\boxtimes$.

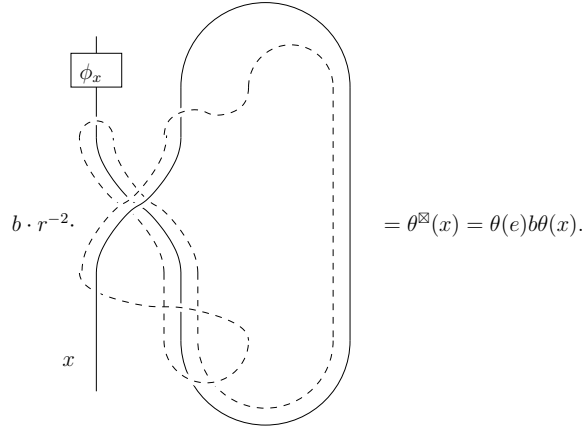
\figref{fig-twist}{The twist of odd $x$ in $\mathcal C^\boxtimes$.}{width=3in}

Twists for even objects have the same value as in $\mathcal C$. The picture for the odd twist is shown in Figure~\ref{fig-twist}. Then

\[\theta^\boxtimes_x = r^{-2} b \frac{\tilde{s}_{e,x}}{\operatorname{x}} \theta_x = -b \theta_x.\]

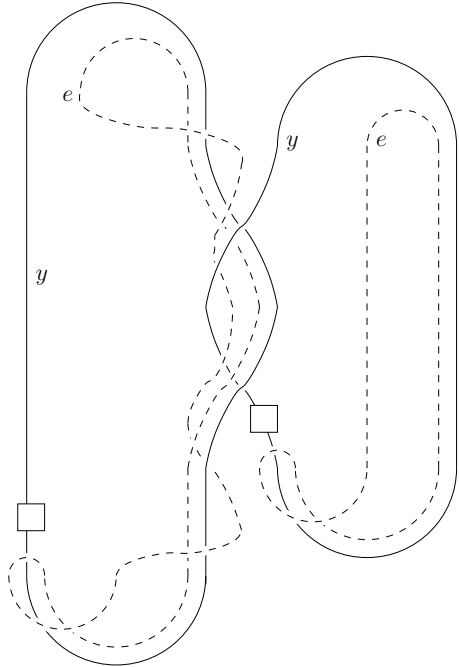
\figref{fig-sentry}{The $S$-matrix entry for odd objects $x$ and $y$.}{height=3.5in}

Let $x$ and $y$ be simple objects in $\mathcal C^\boxtimes$. If either is even, $\tilde{s}^\boxtimes_{x,y} = \tilde{s}_{x,y}$. Otherwise, $\tilde{s}_{x,y}$ is given in Figure~\ref{fig-sentry}. The evaluation is then
\[\tilde{s}^\boxtimes_{x,y} = r^{-4}(-1)^3 b^2 \frac{\tilde{s}_{x^*,e}}{d_x}\tilde{s}_{x,y} = b^2 \tilde{s}_{x,y} = r \tilde{s}_{x,y}.\]

\begin{proposition}
  $(\operatorname{Ob}(\mathcal C), \boxtimes, \beta, \lambda, \rho, c^\boxtimes, \theta^\boxtimes)$ is a modular category.
\end{proposition}
\begin{proof}
  Follows immediately by Theorem \ref{zested modular closure}.
\end{proof}


\end{document}